\documentclass[smallcondensed,numbook]{svjour3}     
\usepackage{graphicx}
\usepackage{amsmath}
\usepackage[nottoc,notlof,notlot]{tocbibind} 
\usepackage[titles]{tocloft}
\usepackage{amssymb}
\usepackage{bbold}
\usepackage{multirow}
\usepackage[utf8]{inputenc}
\usepackage{ushort}
\usepackage{authblk}
\usepackage{leftidx}
\usepackage{amsfonts}
\usepackage{booktabs} 
\usepackage{array} 
\usepackage{paralist} 
\usepackage{verbatim} 
\usepackage{fancyhdr}
\usepackage{slashbox}
\usepackage{here}
\usepackage{anyfontsize}

\usepackage{url}
\newcommand{\vertiii}[1]{{\left\vert\kern-0.25ex\left\vert\kern-0.25ex\left\vert #1 
		\right\vert\kern-0.25ex\right\vert\kern-0.25ex\right\vert}}

\newcommand{\Lnorm}[2]{\left(#1, #2\right)}
\newcommand{\Enorm}[2]{a\left(#1, #2\right)}

\newcommand{\llnorm}[1]{\left\lVert #1\right\rVert_{L_2(\Omega)}}
\newcommand{\edgelnorm}[1]{\left\lVert #1\right\rVert_{L_2(e)}}
\newcommand{\elementlnorm}[1]{\left\lVert #1\right\rVert_{L_2(E)}}
\newcommand{\ilnorm}[1]{\left\lVert #1\right\rVert_{L_\infty(0,T;L_2(\Omega))}}
\newcommand{\Llnorm}[1]{\left\lVert #1\right\rVert_{L_2(0,T;L_2(\Omega))}}
\newcommand{\hnorm}[1]{\left\lVert #1\right\rVert_{H^1(\Omega)}}
\newcommand{\DGnorm}[1]{\left\lVert #1\right\rVert_{\mathcal{V}}}
\newcommand{\enorm}[1]{\left\lVert #1\right\rVert_{V}}

\newcommand{\gnorm}[1]{\left\lVert #1\right\rVert_{L_2(\Gamma_N)}}

\newcommand{\sobolevl}[1]{\vertiii{v}_{H^0(\mathcal{E}_h)}}
\newcommand{\Gnorm}[2]{\left(#1, #2\right)_{L_2(\Gamma_N)}}

\usepackage{amsthm}
\theoremstyle{definition}
\makeatletter
\renewenvironment{proof}[1][\proofname]{\par
	\normalfont
	\topsep6\p@\@plus6\p@ \trivlist
	\item[\hskip5\labelsep\itshape
	#1\@addpunct{. }]
}{%
	\qed\endtrivlist
}
\makeatother
\let\OLDthebibliography\thebibliography
\renewcommand\thebibliography[1]{
	\OLDthebibliography{#1}
	\setlength{\labelsep}{.5em}
}
\newtheorem{eg}{Example}[section]
\newtheorem{rmk}{Remark}[section]

\numberwithin{equation}{section}
\allowdisplaybreaks
\journalname{Advances in Computational Mathematics}
\begin{document}

\title{A priori error analysis for a finite element approximation of dynamic viscoelasticity problems involving a fractional order integro-differential constitutive law\thanks{Jang gratefully acknowledges the supported of a scholarship from Brunel University London. 
}}


\author{Yongseok Jang\textsuperscript{*,a,b}         \and
        Simon Shaw\textsuperscript{a} 
}


\institute{* Corresponding author\\
	\textsuperscript{a} 
	Department of Mathematics, Brunel University London, Uxbridge, UB8 3PH, UK.\\
	\textsuperscript{b}
	\emph{Present address}: CERFACS, 42 Avenue Gaspard Coriolis, 31057, Toulouse, France.\\	
			Y. Jang \at
              \email{yongseok.jang@brunel.ac.uk; yongseok.jang@cerfacs.fr}           
           \and
           S. Shaw \at
              \email{simon.shaw@brunel.ac.uk} 
}

\date{Received: 06 Jul 2020 / Accepted: 04 Mar 2021}

\maketitle

\begin{abstract}
	\noindent We consider a fractional order viscoelasticity problem modelled by a \textit{power-law} type stress relaxation function. This viscoelastic problem is a Volterra integral equation of the second kind with a weakly singular kernel where the convolution integral corresponds to fractional order differentiation/integration. We use a spatial finite element method and a finite difference scheme in time. Due to the weak singularity, fractional order integration in time is managed approximately by linear interpolation so that we can formulate a fully discrete problem. In this paper, we present a stability bound as well as \textit{a priori} error estimates. Furthermore, we carry out numerical experiments with varying regularity of exact solutions at the end.  
\end{abstract}

\section{Introduction}\label{sec:introduction} Materials that exhibit elastic and viscous response are called
viscoelastic materials such as soft tissues, metals at high temperature, and polymers, e.g. see \cite{hunter1976mechanics}. Deformation of a material follows a momentum equation. It is defined by
\begin{equation}
\rho\ddot{\boldsymbol{u}}-\nabla\cdot\ushort{\boldsymbol{\sigma}}=\boldsymbol{f}\qquad\text{on $\Omega\times(0,T]$},\label{eq:primal:visco}
\end{equation}
where $\ddot{\boldsymbol{u}}$ is acceleration, $\nabla\cdot\ushort{\boldsymbol{\sigma}}$ is the divergence of stress, $\boldsymbol{f}$ is an external body force (e.g. see \cite{VE,DGV}), $\Omega$ is a spatial domain in $\mathbb{R}^d$ for $d=1,2,3$ and $(0,T]$ is a time interval domain for $T>0$.  Here we denote first and second time derivative by single and double overdot, respectively, for example, $\dot{\boldsymbol{u}}$ is velocity where we have displacement $\boldsymbol{u}$. A constitutive equation of linear viscoelasticity is formulated as an integro-differential equation which is characterised with a stress relaxation function \cite{hunter1976mechanics,VE,findley2013creep,drozdov1998viscoelastic} such that
\begin{equation}
\ushort{\boldsymbol{\sigma}}(t)=\ushort{\boldsymbol{D}}\varphi(t)\ushort{\boldsymbol{\varepsilon}}(0)+\int^t_0\ushort{\boldsymbol{D}}\varphi(t-s)\ushort{\dot{\boldsymbol{\varepsilon}}}(s)ds,\label{eq:consti:visco}
\end{equation}
where $\ushort{\boldsymbol{D}}$ is a fourth order symmetric positive definite tensor, for example 
\[D_{ijkl}=D_{jikl}=D_{ijlk}=D_{klij},\] $\ushort{\boldsymbol{\varepsilon}}$ is the strain, and the form of $\varphi$ depends on which viscoelastic model is invoked. A rheological models such as the Maxwell, Voigt and Zener models exhibit exponentially decaying stress relaxation \cite{VE}. For more details, see \cite{findley2013creep,drozdov1998viscoelastic,golden2013boundary} and references therein. In the case of generalised Maxwell model, quasi-static and dynamic linear viscoelastic problems have been dealt with by finite element approximation in  \cite{DGV,riviere2003discontinuous,shaw1998numerical,shaw1999numerical,jang2020finite}.

Another choice of a stress relaxation function, namely \textit{power-law}, was employed by Nutting \cite{nutting1921new}, e.g. see also \cite{torvik1984appearance,koeller1984applications}. The {power-law} type kernel has been naturally introduced in an intermediate sense between elasticity and viscosity \cite{VE}. To be specific, classical continuum mechanics provides that $\ushort{\boldsymbol{\sigma}}\propto\ushort{\boldsymbol{\varepsilon}}$ in elastic solid and $\ushort{\boldsymbol{\sigma}}\propto\ushort{\dot{\boldsymbol{\varepsilon}}}$ in viscous liquid so that the constitutive relation of viscoelasticity could exhibit $\ushort{\boldsymbol{\sigma}}\propto\partial_t^{\alpha}\ushort{{\boldsymbol{\varepsilon}}}$, where $\partial_t^{\alpha}$ is the fractional $\alpha$ order differential operator such that \[\partial_t^{\alpha}\ushort{{\boldsymbol{\varepsilon}}}(t)=\frac{1}{\Gamma(1-\alpha)}\frac{\partial}{\partial t}\int^t_0(t-s)^{-\alpha}\ushort{{\boldsymbol{\varepsilon}}}(s)ds,\] for $\alpha\in(0,1)$ and $\Gamma$ is the Gamma function. For example, it is argued that $\ushort{\boldsymbol{\sigma}}\propto\partial_t^{0.56}\ushort{{\boldsymbol{\varepsilon}}}$ in elastomer 3M-467 in \cite{torvik1984appearance}. In this manner, the power-law type stress relaxation kernel for viscoelasticity \cite{VE} is introduced by
	\begin{equation}
	\varphi(t)=\varphi_0+\varphi_1t^{-\alpha} \label{eq:stressrelax:power},
	\end{equation}where $\varphi_0$ is non-negative, $\varphi_1$ is positive and $\alpha\in(0,1)$. Consequently, the power-law type kernel leads us to derive a fractional order viscoelastic model. 

 Due to the weakly singular kernel in the fractional order viscoelasticity model, the standard quadrature rules such as the trapezoidal rule, are unable to work. For instance, the typical quadrature rules require function values of the integrand for all nodes but $\varphi(0)$ is unbounded.
Hence we need to find alternative methods which resolve the singularity at $t=0$. In \cite{li2011numerical}, some numerical approach for fractional calculus was introduced based on interpolation techniques with various accuracy orders.
McLean and Thom\'ee \cite{mclean1993numerical,mclean2010numerical,mclean2010maximum} developed numerical analysis of a fractional order evolution equation which is a scalar analogue of a fractional order viscoelasticity problem of power-law type, and they presented error analysis with the homogeneous Dirichlet boundary condition. 

In this article, we study the fractional order viscoelastic model problem with mixed boundary conditions. We consider finite element approximation for the fractional order viscoelastic model given by power-law type stress relaxation. On account of the weak singularity, we may encounter some difficulty in \textit{a priori} analysis. To resolve this issue, we introduce the linear interpolation technique \cite{li2011numerical,linz2001theoretical} while we employ spatial finite element method and Crank-Nicolson finite difference method in time. We show stability bounds as well as spatially optimal error bounds but without Gr\"onwall inequality for time integral not to produce exponentially increasing bounds in time. In terms of the weak singularity, we will discuss regularity of solutions to obtain suboptimal and optimal convergence orders with respect to time.

Here, we would like to highlight that the well-posedness for the fractional order integro-differential equation with the mixed boundary condition can be shown by introducing Markov's inequality but without Gr\"onwall's inequality. Despite the weak singularity in the power-law type kernel and limitations for higher regularity of solutions in time, the fully discrete solutions have better order of accuracy than first order schemes. We can prove it by means of duality arguments and $L_\infty$ approach in time rather than by the use of Gr\"onwall's inequality and spectral method.

This article is arranged as follows. In Section \ref{sec:preliminary}, we introduce fundamental definitions of fractional calculus, the finite element method and our notation. In Section \ref{sec:modelproblem}, we give more suppositions to derive the reduced model of \eqref{eq:primal:visco} and define discrete formulations. In Section \ref{sec:analysis}, we state and prove stability bounds as well as error estimates. By using the fully discrete formula, numerical experiments are carried out in FEniCS Project (\texttt{https://fenicsproject.org/}) in Section \ref{sec:numeric}. At the end, we conclude with Section \ref{sec:conclusion}.
\section{Preliminary}\label{sec:preliminary}
\noindent
According to \cite{oldham1974fractional,malinowska2012introduction,miller1993introduction}, we define 
	Riemann-Liouville fractional integral as follows.
	If $f\in L_1[a,b]$, the $\alpha$ order integral of $f$ is given by
	\[\leftidx{_{a}}I_t^\alpha f(t)=\frac{1}{\Gamma(\alpha)}\int^t_af(s)(t-s)^{\alpha-1}ds,\ t>a,\]where $\alpha$ is positive. We can also rewrite the fractional integral in convolution form as
	\[\leftidx{_{a}}I_t^\alpha f(t)=\beta_\alpha*f(t),\]
where $\beta_\alpha(t)={t^{\alpha-1}}/{\Gamma(\alpha)}$ and $*$ denotes Laplace convolution such that \[f_1*f_2(t)=\int^{t}_{0}f_1(t-s)f_2(s)ds.\] Note that $\beta_\alpha(t)$ is a weakly singular kernel for $0<\alpha<1$.

We introduce and use some standard notations so that the usual $L_p(\Omega),\ H^s(\Omega)$ and $W^s_p(\Omega)$ denote Lebesgue, Hilbert and Sobolev space, respectively, where $s$ and $p$ are non-negative. For any Banach space $X$, $\lVert\cdot\rVert_X$ is the $X$ norm, for example, $\llnorm{\cdot}$ is the $L_2(\Omega)$ norm induced by the $L_2(\Omega)$ inner product which we denote for brevity by $(\cdot,\cdot)$, but for $S\subset\bar\Omega$, $(\cdot,\cdot)_{L_2(S)}$ is the $L_2(S)$ inner product. In case of time dependent functions, we expand this notation such that if $f\in L_p(0,T;X)$ for some Banach space $X$, we define
\[\lVert f\rVert_{L_p(0,t_0;X)}=\left(\int^{t_0}_0\lVert f(t)\rVert_X^pdt\right)^{1/p}\]
for $t_0\leq T$ and $1\leq p<\infty$. When $p=\infty$, we shall use \textit{essential supremum} norm where $$\lVert f\rVert_{L_\infty(0,t_0;X)}=\mathop{\mathrm{ess~sup}}\limits_{0\leq t\leq t_0}\lVert f(t)\rVert_X .$$ Also, we define \textit{H\"older} norm for $f\in C^s(0,T;X)$ by $$\lVert f\rVert_{C^s(0,T;X)}=\mathop{\mathrm{max}}\limits_{0\leq k\leq s}\mathop{\mathrm{sup}}\limits_{0\leq t\leq T}\left\lVert \frac{\partial^k}{\partial t^k}f(t)\right\rVert_X.$$

Let us define a framework for our finite element method. We assume that $\Omega$ is an open bounded convex polytopic domain, $\Gamma_D$ is the positive measured \textit{Dirichlet} boundary, and the \textit{Neumann} boundary $\Gamma_N$ is given by $\Gamma_N=\partial\Omega\backslash\Gamma_D$. 
For use later we recall the trace inequality, \begin{align}
\lVert v\rVert_{L_2(\partial\Omega)}&\leq C\hnorm {v}, \textrm{ for any }v\in H^1(\Omega),
\label{cgtrace}
\end{align}where $C$ is a positive constant depending only on $\Omega$ and its boundary.

Let $V$ be a subspace of $H^1(\Omega)$ such that \begin{align*}
V=\left\{v\in H^1(\Omega)\ |\ v=0\ \textrm{ on } {\Gamma_D}\right\},
\end{align*}and $V^h$ be a finite element space of polynomial of degree $k$ in $V$. In particular, we consider conforming meshes and Lagrange finite elements for the construction of the finite element space \cite{MTE,wheeler}.
For the sake of our model problem, we will use similar notations for vector-valued functions. Let us define \begin{align*}
\boldsymbol{V}=[V]^d=\left\{\boldsymbol{v}\in [H^1(\Omega)]^d\ |\ \boldsymbol v=\boldsymbol 0\ \textrm{ on } {\Gamma_D}\right\},
\end{align*}and $\boldsymbol{V}^h=[V^h]^d$. Also, we use inner products of vector-valued (tensor-valued) functions with same notations as scalar cases. For instance, we have
\[\Lnorm{\boldsymbol{v}}{\boldsymbol{w}}=\int_\Omega \boldsymbol{v}\cdot\boldsymbol{w}\ d\Omega,\qquad\Lnorm{\ushort{\boldsymbol{\psi}}}{\ushort{\boldsymbol{\zeta}}}=\int_\Omega \ushort{\boldsymbol{\psi}}:\ushort{\boldsymbol{\zeta}}\ d\Omega,\]for vector-valued functions $\boldsymbol{v}$ and $\boldsymbol{w}$, and second order tensors $\ushort{\boldsymbol{\psi}}$ and $\ushort{\boldsymbol{\zeta}}$.
\section{Model Problem}\label{sec:modelproblem}
Consider the viscoelasticity model problem with the power-law type constitutive relation. Then we have
\begin{eqnarray}
\rho\ddot{\boldsymbol{u}}(t)-\nabla\cdot\left(\ushort{\boldsymbol{D}}(\varphi_0+\varphi_1t^{-\alpha})\ushort{\boldsymbol{\varepsilon}}(\boldsymbol{u}(0) )+\int^t_0\ushort{\boldsymbol{D}}(\varphi_0+\varphi_1(t-s)^{-\alpha})\ushort{\dot{\boldsymbol{\varepsilon}}}(\boldsymbol{u}(s))ds\right)=\boldsymbol{f}(t),\label{eq:powerlaw:visco:eq1}
\end{eqnarray}where $t\in(0,T],\ \alpha\in(0,1)$ and $\ushort{\boldsymbol{\varepsilon}}$ is Cauchy infinitesimal tensor defined by, for any $\boldsymbol{v}\in [H^1(\Omega)]^d$,
\[\varepsilon_{i,j}(\boldsymbol{v})=\frac{1}{2}\left(\frac{\partial v_i}{\partial x_j}+\frac{\partial v_j}{\partial x_i}\right),\qquad\text{for }i,j=1,\ldots,d.\]
Note that the strain tensor $\ushort{\boldsymbol{\varepsilon}}$ is a symmetric second order tensor. Hence we have
\begin{align}
\ushort{\boldsymbol \sigma}:\ushort{\boldsymbol \varepsilon}(\boldsymbol v)=\ushort{\boldsymbol \sigma}:\nabla\boldsymbol v,\qquad\forall\boldsymbol{v}\in [H^1(\Omega)]^d,\label{eq:sym:div}
\end{align}since the stress and strain are symmetric. Using the fundamental theorem of calculus, we can observe that
\[\ushort{\boldsymbol{D}}\varphi_0\ushort{\boldsymbol{\varepsilon}}(\boldsymbol{u}(0) )+\int^t_0\ushort{\boldsymbol{D}}\varphi_0\ushort{\dot{\boldsymbol{\varepsilon}}}(\boldsymbol{u}(s))ds=\ushort{\boldsymbol{D}}\varphi_0\ushort{\boldsymbol{\varepsilon}}(\boldsymbol{u}(t) ),\] and this is purely elastic response.
To simplify \eqref{eq:powerlaw:visco:eq1}, we assume that $\boldsymbol{u}(0)=\boldsymbol{0}$ and $\varphi_0=0$. For convenience of notation, we also define $\hat{\ushort{\boldsymbol{D}}}$ such that
\[\hat{\ushort{\boldsymbol{D}}}=\varphi_1\Gamma(1-\alpha)\ushort{\boldsymbol{D}}.\] Once we denote the velocity vector by $\boldsymbol{w}=\dot{\boldsymbol{u}}$, we can reduce \eqref{eq:powerlaw:visco:eq1} to a lower order differential problem by using fractional integral notation. Thus, we will consider the following model problem: find $\boldsymbol{w}$ such that
\begin{align}
\rho\dot{\boldsymbol{w}}(t)-\nabla\cdot\leftidx{_0}I^{1-\alpha}_t({\hat{\ushort{\boldsymbol{D}}}}\ushort{{\boldsymbol{\varepsilon}}}(\boldsymbol{w}(t)))&=\boldsymbol{f}(t),&\textrm{ on }(0,T]\times\Omega,\label{eq:frac:primal:eq1}\\
\leftidx{_0}I^{1-\alpha}_t({\hat{\ushort{\boldsymbol{D}}}}\ushort{{\boldsymbol{\varepsilon}}}(\boldsymbol{w}(t)))\cdot\boldsymbol{n}&=\boldsymbol{g}_N(t),&\textrm{ on }[0,T]\times\Gamma_N,\label{eq:frac:primal:neumann}\\
\boldsymbol{w}(t)&=\boldsymbol{0},&\textrm{ on }[0,T]\times\Gamma_D,\label{eq:frac:primal:dirichlet}\\
\boldsymbol{w}(0)&=\boldsymbol{w}_0,&\textrm{ on }\Omega,\label{eq:frac:primal:initial}
\end{align}where $\alpha\in(0,1)$, $\hat{\ushort{\boldsymbol{D}}}$ is a symmetric positive definite piecewise constant fourth order tensor and $\boldsymbol{n}$ is an outward unit normal vector. In continuum mechanics, $\boldsymbol{g}_N$ is called \textit{traction}, which is equivalent to $\ushort{\boldsymbol{\sigma}}\cdot\boldsymbol{n}$.
\subsection{Weak Formulation}\label{subsec:semi}
\noindent As taking into account multiplying $\boldsymbol{v}\in\boldsymbol{V}$ by \eqref{eq:frac:primal:eq1} and integrating it over $\Omega$, we are able to obtain the following weak problem: find a mapping $\boldsymbol{w}:[0,T]\mapsto\boldsymbol{V}$ such that 
\begin{gather}
\Lnorm{\rho\dot{\boldsymbol{{w}}}(t)}{\boldsymbol{v}}+\Enorm{\leftidx{_0}I^{1-\alpha}_t\boldsymbol{w}(t)}{\boldsymbol{v}}=F(t;\boldsymbol{v}),\ \forall t\in(0,T],\label{eq:frac:T:primal:eq1}\\
\Enorm{\boldsymbol{w}(0)}{\boldsymbol v}=\Enorm{\boldsymbol{w}_0}{\boldsymbol v},\label{eq:frac:T:primal:eq2}
\end{gather} for any $\boldsymbol{v}\in\boldsymbol{V}$ where $\Enorm{\cdot}{\cdot}$ and $F$ are defined by 
\[
\Enorm{\boldsymbol{w}}{\boldsymbol v}=\int_\Omega{\hat{\ushort{\boldsymbol{D}}}\ushort{\boldsymbol{\varepsilon}}(\boldsymbol{w})}:{\ushort{\boldsymbol{\varepsilon}}(\boldsymbol{v})}d\Omega\]and
\[F(t;\boldsymbol v)=\Lnorm{\boldsymbol{f}(t)}{\boldsymbol{v}}+\Gnorm{\boldsymbol{g}_N(t)}{\boldsymbol{v}}.\]
It is easily to show that \eqref{eq:frac:T:primal:eq1} is a weak form of \eqref{eq:frac:primal:eq1} by \eqref{eq:sym:div} and integration by parts. Straightforwardly, \eqref{eq:frac:primal:initial} gives \eqref{eq:frac:T:primal:eq2}, since the bilinear form is well-defined.

\begin{rmk}
	As is usual in variational problems, we may want to show continuity and coercivity of the bilinear form, and continuity of the linear form. According to Korn's inequality \cite{MTE,ciarlet2010korn,horgan1983inequalities,nitsche1981korn}, 
	\[C\hnorm{\boldsymbol v}^2\leq \lambda_{\textrm{min}}\llnorm{\ushort{\boldsymbol{\varepsilon}}(\boldsymbol{v})}^2\leq\Enorm{\boldsymbol{v}}{\boldsymbol{v}}\leq \lambda_{\textrm{max}}\llnorm{\ushort{\boldsymbol{\varepsilon}}(\boldsymbol{v})}^2\leq \lambda_{\textrm{max}}\hnorm{\boldsymbol{v}}^2\]for any $\boldsymbol{v}\in \boldsymbol{V}$ where $C$ is a positive constant independent of $\boldsymbol{v}$ and ($\lambda_{\textrm{min}},\lambda_{\textrm{max}}$) is a pair of the minimum and maximum eigenvalues of $\hat{\ushort{\boldsymbol{D}}}$. Also, we can observe that
	\[|\Enorm{\boldsymbol{v}}{\boldsymbol{w}}|\leq\lambda_{\textrm{max}}\llnorm{\ushort{\boldsymbol{\varepsilon}}(\boldsymbol{v})}\llnorm{\ushort{\boldsymbol{\varepsilon}}(\boldsymbol{w})}\nonumber\\
	\leq \lambda_{\textrm{max}}\hnorm{\boldsymbol{v}}\hnorm{\boldsymbol{w}},\]
	for any $\boldsymbol{v},\boldsymbol{w}\in\boldsymbol{V}$. Therefore, the bilinear form is coercive and continuous. 
	Furthermore, when we define the energy norm $\enorm{\cdot}$ on $\boldsymbol{V}$ by
	\[\enorm{\boldsymbol{v}}=\sqrt{\Enorm{\boldsymbol{v}}{\boldsymbol{v}}},\]
	we can observe norm equivalence between the $H^1$ and energy norms on $\boldsymbol{V}$ and we have\begin{align}
	|\Enorm{\boldsymbol{v}}{\boldsymbol{w}}|\leq\enorm{\boldsymbol{v}}\enorm{\boldsymbol{w}},\label{eq:conti}
	\end{align}for any $\boldsymbol{v},\boldsymbol{w}\in\boldsymbol{V}$. On the other hand, the use of Cauchy-Schwarz inequality and trace inequality allows us to show the continuity of the linear form.
\end{rmk}
\begin{rmk}According to \cite{li2011developing}, $t^{-\alpha}$ for $0<\alpha<1$ is a positive definite kernel such that for $T>0$
	\begin{equation}
	\int^T_0\phi(t)\int^t_0(t-s)^{-\alpha}\phi(s)dsdt=\int^T_0\int^t_0(t-s)^{-\alpha}\phi(s)\phi(t)dsdt\geq 0,\ \forall\phi\in C[0,T],\label{eq:frac:kernel:pd}
	\end{equation}and hence
	\begin{equation}
	\int^T_0\leftidx{_0}I^{1-\alpha}_t\phi(t)\phi(t)dt=\frac{1}{\Gamma(1-\alpha)}\int^T_0\int^t_0(t-s)^{-\alpha}\phi(s)\phi(t)dsdt\geq 0.\label{eq:frac:int:pd}
	\end{equation}
\end{rmk}
In order to carry out stability analysis, we shall use \eqref{eq:frac:int:pd}.
\begin{theorem}\label{thm:frac:cg:semi:sta}
	Suppose that $\boldsymbol{w}\in L_\infty(0,T;\boldsymbol{V})\cap H^1(0,T;[L_2(\Omega)]^d)$, $\boldsymbol{f}\in L_2(0,T;[L_2(\Omega)]^d)$ and $\boldsymbol{w}_0\in\boldsymbol{V}$. In addition, we assume $\boldsymbol{g}_N=\boldsymbol{0}$. Then there exists a positive constant $C$ such that 
	\begin{align*}
	{\rho}\ilnorm{\boldsymbol{w}}^2	\leq& CT\left({\rho}\enorm{\boldsymbol{w}_0}^2+\Llnorm{\boldsymbol{f}}^2\right).
	\end{align*}
\end{theorem}
\begin{proof}
	Let $\boldsymbol{v}=\boldsymbol{w}(t)$ in \eqref{eq:frac:T:primal:eq1} to get
	\begin{align}
	\frac{\rho}{2}\frac{d}{dt}\llnorm{\boldsymbol{w}(t)}^2+\Enorm{\leftidx{_0}I^{1-\alpha}_t\boldsymbol{w}(t)}{\boldsymbol{w}(t)}=F(\boldsymbol{w}(t)).\label{eq:frac:cg:pf:semi:sta:eq1}
	\end{align}
	Taking into account the second term of the left hand side of \eqref{eq:frac:cg:pf:semi:sta:eq1}, the definition of the fractional integral gives
	\begin{align}
	\Enorm{\leftidx{_0}I^{1-\alpha}_t\boldsymbol{w}(t)}{\boldsymbol{w}(t)}=\frac{1}{\Gamma(1-\alpha)}\int^t_0(t-t')^{-\alpha}\Enorm{\boldsymbol{w}(t')}{\boldsymbol{w}(t)}dt',
	\label{eq:frac:cg:pf:semi:sta:eq2}\end{align}by Leibniz integral rule.
	By substitution of \eqref{eq:frac:cg:pf:semi:sta:eq2} into \eqref{eq:frac:cg:pf:semi:sta:eq1}, integrating over time yields
	\begin{align}
	\frac{\rho}{2}&(\llnorm{\boldsymbol{w}(\tau)}^2-\llnorm{\boldsymbol{w}(0)}^2)+\frac{1}{\Gamma(1-\alpha)}\int^\tau_0\int^t_0(t-t')^{-\alpha}\Enorm{\boldsymbol{w}(t')}{\boldsymbol{w}(t)}dt'dt\nonumber\\
	=&\int^\tau_0F(\boldsymbol{w}(t))dt,\label{eq:frac:cg:pf:semi:sta:eq3}
	\end{align}for $0<\tau\leq T$.
	In the double integral, we can expand the bilinear form and take spatial integration outside so that \eqref{eq:frac:kernel:pd} gives
	\begin{align*}
	\int_0^\tau\int^t_0(t-t')^{-\alpha}\Enorm{\boldsymbol{w}(t')}{\boldsymbol{w}(t)}dt'dt
	=\int_\Omega\int_0^\tau\int^t_0(t-t')^{-\alpha}\hat{\ushort{\boldsymbol{D}}}^{1/2}\ushort{{\boldsymbol{\varepsilon}}}({\boldsymbol{w}(t')}):\hat{\ushort{\boldsymbol{D}}}^{1/2}\ushort{{\boldsymbol{\varepsilon}}}({\boldsymbol{w}(t)}) dt'dtd\Omega
	\geq 0,
	\end{align*}where $\hat{\ushort{\boldsymbol{D}}}^{1/2}$ is a symmetric positive definite fourth order tensor satisfying $\hat{\ushort{\boldsymbol{D}}}^{1/2}\hat{\ushort{\boldsymbol{D}}}^{1/2}=\hat{\ushort{\boldsymbol{D}}}$ by the use of spectral decomposition. As a consequence \eqref{eq:frac:cg:pf:semi:sta:eq3} yields
	\begin{align}
	\frac{\rho}{2}\llnorm{\boldsymbol{w}(\tau)}^2
	\leq&\frac{\rho}{2}\llnorm{\boldsymbol{w}(0)}^2+\int^\tau_0F(\boldsymbol{w}(t))dt.\label{eq:frac:cg:pf:semi:sta:eq4}
	\end{align}
We can observe a bound of the last term	in \eqref{eq:frac:cg:pf:semi:sta:eq4} such that
	\begin{align*}
	\int^\tau_0F(\boldsymbol{w}(t))dt\leq&\int^\tau_0\llnorm{\boldsymbol{f}(t)}\llnorm{\boldsymbol{w}(t)}dt\\
	\leq&\frac{\epsilon_a}{2}\ilnorm{\boldsymbol{w}}^2+\frac{T}{2\epsilon_a}\Llnorm{\boldsymbol{f}}^2
	\end{align*}by Cauchy-Schwarz and Young's inequalities for any positive $\epsilon_a$. Since $\tau$ is arbitrary, we can complete the proof by choice of $\epsilon_a=\rho/2$ and therefore we have\begin{align}
	\frac{\rho}{4}\ilnorm{\boldsymbol{w}}^2
	\leq&CT\bigg(\rho\llnorm{\boldsymbol{w}(0)}^2+\Llnorm{\boldsymbol{f}}^2\bigg),\label{eq:frac:cg:pf:semi:sta:eq5}
	\end{align}where $C$ is a positive constant. Moreover, it is seen that $\llnorm{\boldsymbol{w}(0)}^2\leq C\enorm{\boldsymbol{w}_0}^2$ by coercivity and norm equivalence in \eqref{eq:frac:T:primal:eq2},
	hence the theorem is proved.
\end{proof}
Since $\boldsymbol{V}^h\subset\boldsymbol{V}$ is a finite dimensional subspace, Theorem \ref{thm:frac:cg:semi:sta} holds for $\boldsymbol{V}^h$. It means we can find a semidiscrete solution $[0,T]\mapsto\boldsymbol{V}^h$ which fulfils \eqref{eq:frac:T:primal:eq1}-\eqref{eq:frac:T:primal:eq2} with the stability bound in Theorem \ref{thm:frac:cg:semi:sta}. 
\subsection{Fully Discrete Formulation}\label{subsec:fully}
\noindent Next, we are going to formulate a fully discrete problem. We use the Crank-Nicolson finite difference scheme for time discretization but it is also necessary to introduce numerical methods for fractional order integral.

Let $\Delta t=T/N$ for some $N\in\mathbb{N}$. Define $t_n=n\Delta t$ for $n=0,\ldots,N$ and denote our fully discrete solution by $\boldsymbol{W}_h^{n}\in\boldsymbol{V}^h$ for $n=0,\ldots,N$. In the way of the Crank-Nicolson finite difference method, we will approximate first time derivatives by
\[\frac{\dot{\boldsymbol{w}}(t_{n+1})+\dot{\boldsymbol{w}}(t_{n})}{2}\approx\frac{\boldsymbol{W}_h^{n+1}-\boldsymbol{W}_h^{n}}{\Delta t}\qquad\text{for $n=0,\ldots,N-1$.}\]
Due to the weak singularity in the fractional integral, we should be cautious when using numerical integration. We will use linear interpolation technique from \cite{li2011numerical}, and define the piecewise linear interpolation of $\boldsymbol{w}$ such that for $n=1,\ldots,N$,
\[\bar{\boldsymbol{w}}_n(t)=-\frac{t-t_{n}}{\Delta t}\boldsymbol{w}(t_{n-1})+\frac{t-t_{n-1}}{\Delta t}\boldsymbol{w}(t_n)\text{ where }t\in[t_{n-1},t_n].\]If $\boldsymbol{{w}}$ is of $C^2$ in time, we have for $t\in[t_{n-1},t_n]$,
\[\boldsymbol{E}_n(t):=\boldsymbol{w}(t)-\bar{\boldsymbol{w}}_n(t)=\frac{1}{2}\ddot{\boldsymbol w}(\xi_t)(t-t_{n-1})(t-t_{n})\text{ for some }\xi_t\in[t_{n-1},t_n],\]
by Rolle's theorem. If $\boldsymbol{w}(t)\in[H^s(\Omega)]^d$ for any $t\in[t_{n-1},t_n]$, it holds that
\begin{align}\lVert\boldsymbol{E}_n(t)\rVert_{H^s(\Omega)}\leq\frac{\Delta t^2}{2}\lVert\ddot{\boldsymbol w}\rVert_{C^0(t_{n-1},t_n;H^s(\Omega))}\label{eq:num:linear:conti}\end{align}
Then we can obtain the following numerical approximation 
\begin{align}\leftidx{_0}I^{1-\alpha}_{t}\boldsymbol{w}(t_n)=&\frac{1}{\Gamma(1-\alpha)}\sum_{i=1}^{n}\int^{t_i}_{t_{i-1}}\left(\bar{\boldsymbol{w}}_i(t')+\boldsymbol{E}_i(t')\right)(t_n-t')^{-\alpha}dt'\nonumber\\=&\frac{\Delta t^{1-\alpha}}{\Gamma(3-\alpha)}\sum_{i=0}^{n}B_{n,i}\boldsymbol{w}(t_i)+\frac{1}{\Gamma(1-\alpha)}\sum_{i=1}^{n}\int^{t_i}_{t_{i-1}}\boldsymbol{E}_i(t')(t_n-t')^{-\alpha}dt'\nonumber\\:=&\boldsymbol{q}_n(\boldsymbol{w})+\frac{1}{\Gamma(1-\alpha)}\sum_{i=1}^{n}\int^{t_i}_{t_{i-1}}\boldsymbol{E}_i(t')(t_n-t')^{-\alpha}dt',\label{eq:numeric_int}\end{align} where
\begin{align*}
B_{n,i}=\left\{ \begin{array}{cc}
n^{1-\alpha}({2-\alpha}-n)+(n-1)^{{2-\alpha}},&i=0,\\
(n-i-1)^{{2-\alpha}}+(n-i+1)^{{2-\alpha}}-2(n-i)^{{2-\alpha}},&i=1,\ldots,n-1,\\
1,&i=n.
\end{array}\right.
\end{align*}Note that $0<B_{n,i}<2$ for any $n$ and $i=0,\ldots,n$. By using Cauchy-Schwarz inequality, if $\boldsymbol{w}\in C^2(0,T;[H^s(\Omega)]^d)$, we can derive the numerical error such that by \eqref{eq:num:linear:conti},
\begin{align}
\lVert\leftidx{_0}I^{1-\alpha}_{t}\boldsymbol{w}(t_n)-\boldsymbol{q}_n(\boldsymbol{w})\rVert_{H^s(\Omega)}\leq& \frac{\Delta t^2}{2\Gamma(1-\alpha)}\lVert\boldsymbol{\ddot{\boldsymbol w}}\rVert_{C^0(0,t_n;H^s(\Omega))}\int^{t_n}_0(t_n-t')^{-\alpha}dt'\nonumber\\
\leq& \frac{T^{1-\alpha}}{2\Gamma(2-\alpha)}\lVert\boldsymbol{\ddot{\boldsymbol w}}\rVert_{C^0(0,T;H^s(\Omega))}\Delta t^2\label{eq:linearinterpolation:error}.
\end{align}

Consequently, the use of Crank-Nicolson method and the numerical integration leads us to obtain the fully discrete formulation as follows:
find $\boldsymbol{W}_h^n\in \boldsymbol{V}^h$ for $n=0,\ldots,N$ such that for any $\boldsymbol{v}\in \boldsymbol{V}^h$,
\begin{gather}
\Lnorm{\rho\frac{\boldsymbol{W}_h^{n+1}-\boldsymbol{W}_h^n}{\Delta t}}{\boldsymbol{v}}+\Enorm{\frac{\boldsymbol{q}_{n+1}(\boldsymbol{W}_h)+\boldsymbol{q}_n(\boldsymbol{W}_h)}{2}}{\boldsymbol{v}}=\frac{1}{2}(F(t_{n+1};\boldsymbol{v})+F(t_n;\boldsymbol{v})),\label{eq:frac:T:fully:eq1}
\end{gather}$\forall n=0,\ldots,N-1$, and
\begin{gather}
\Enorm{\boldsymbol{W}_h^0}{\boldsymbol v}=\Enorm{\boldsymbol{w}_0}{\boldsymbol v}.\label{eq:frac:T:fully:eq2}
\end{gather}
\section{Stability and Error Analysis}\label{sec:analysis}
\noindent As shown in Theorem \ref{thm:frac:cg:semi:sta}, we will carry out stability analysis in a fully discrete sense. One can show a stable bound then the existence and uniqueness of the discrete solution is possessed simultaneously. In a similar way with the stability analysis, we can derive error estimates by introducing the elliptic projection.
\subsection{A Stability Bound}\label{subsec:stability}
\noindent Let us present the following inverse polynomial trace theorem and Markov inequality.
\begin{theorem}\label{thm:inversetrace}
	\textnormal{\textbf{Inverse Polynomial Trace Theorem} \cite{WARBURTON20032765}}\\
	Let $E$ be a triangle in 2D or a tetrahedron in 3D and $e$ be an edge in 2D or a face in 3D of $E$. Suppose $\mathcal{P}_k(E)$ is a set of polynomials of degree $k$ on $E$. Then there exists a trace inequality such that 
	\begin{align*}
	\forall v\in\mathcal{P}_k(E),\ \forall e\subset\partial E,\ &\edgelnorm{v}\leq Ch_E^{-1/2}\elementlnorm{v},
	\end{align*}where $h_E$ is a diameter of $E$ and $C$ is a positive constant and is independent of $h_E$ but depending on the degree of polynomial $k$ and the dimension $d$.
\end{theorem}\vspace{0.5cm}

\begin{theorem}\label{thm:inverseinequality}
	\textnormal{\textbf{Inverse Inequality(or Markov Inequality)} \cite{DG,ozisik2010constants}}\\
	For any element $E$, there is a positive constant $C$ such that 
	for any $0\leq j\leq k$,
	\begin{align*}
	\forall v\in \mathcal{P}_k(E),\ \elementlnorm{\nabla^jv}\leq Ch_E^{-j}\elementlnorm{v},\ \text{where }\nabla^j=\left\{
	\begin{array}{cc}
	\nabla(\nabla^{j-1}),&\text{if $j$ odd},\\
	\nabla\cdot(\nabla^{j-1}),&\text{if $j$ even}.
	\end{array}\right. 
	\end{align*}
\end{theorem}\vspace{0.5cm}

If we assume a quasi-uniform mesh, then we can derive 
\begin{align}
\gnorm{\boldsymbol v}^2\leq Ch^{-1}\llnorm{\boldsymbol v}^2,\qquad
\forall \boldsymbol{v}\in\boldsymbol{V}^h,\label{eq:trace:inv} 
\end{align}and
\begin{align}
\llnorm{\ushort{\boldsymbol\varepsilon}(\boldsymbol v)}\leq Ch^{-1}\llnorm{\boldsymbol v},\qquad
\forall \boldsymbol{v}\in\boldsymbol{V}^h.\label{eq:inv} 
\end{align}
Usually to estimate trace terms, trace inequalities, e.g. \eqref{cgtrace}, are used rather than the inverse polynomial trace theorem \cite{MTE}. The typical trace inequality contains the $H^1(\Omega)$ norm, however our problem has a difficulty in dealing with $H^1(\Omega)$ norm due to numerical integration of fractional order integral $\boldsymbol{q}_n$. To be specific, since we can only derive the energy norm of the fractional integrals in stability analysis, we are unable to manage the trace norm of the discrete solution, whereas \eqref{eq:trace:inv} allows us to analyse the trace terms in $L_2(\Omega)$ norm sense. Moreover, we note that the inverse polynomial trace theorem can be employed only in polynomial spaces, which means that \eqref{eq:trace:inv} does not hold in $\boldsymbol{V}$ so we supposed $\boldsymbol{g}_N=\boldsymbol{0}$ in Theorem \ref{thm:frac:cg:semi:sta}. Hereafter, we assume $\boldsymbol{V}^h$ is constructed with a quasi-uniform mesh and hence we can deal with non-zero $\boldsymbol{g}$ as well.

\begin{theorem}\label{thm:frac:cg:fully:sta}
	Suppose $\boldsymbol{f}\in L_2(0,T;[L_2(\Omega)]^d),\ \boldsymbol{ g}_N\in L_2(0,T;[L_2(\Gamma_N)]^d), \textrm{ and }\boldsymbol{w}_0\in\boldsymbol{V}$. Then there exists a unique discrete solution of \eqref{eq:frac:T:fully:eq1} and \eqref{eq:frac:T:fully:eq2} such that 
		\begin{align*}
	&\max_{0\leq n\leq N}\llnorm{\boldsymbol{W}_h^{n}}^2+\frac{\Delta t^{2-\alpha}}{\Gamma(3-\alpha)}\sum_{n=0}^{N-1}\enorm{\boldsymbol{W}_h^{n+1}+\boldsymbol{W}_h^n}^2
	\nonumber\\\ &\leq CT\bigg(\enorm{\boldsymbol{w}_{0}}^2+\Delta t\sum\limits_{n=0}^{N}\llnorm{\boldsymbol{f}(t_n)}^2+\Delta t\sum\limits_{n=0}^{N}h^{-1}\gnorm{{\boldsymbol{g}_N(t_n)}}^2\bigg),
	\end{align*}where $C$ is independent of the solution, $\Delta t$ and $h$.
\end{theorem}
\begin{proof}
	Let $m\in\{1,\ldots, N\}$. A choice of $\boldsymbol{v}=2\Delta t(\boldsymbol{W}_h^{n+1}+\boldsymbol{W}_h^n)$ in \eqref{eq:frac:T:fully:eq1} and summation from $n=0$ to $n=m-1$ yields	
	\begin{align}
	2\rho&\left(\llnorm{\boldsymbol{W}_h^{m}}^2-\llnorm{\boldsymbol{W}_h^{0}}^2\right)+\Delta t\sum_{n=0}^{m-1}\Enorm{\boldsymbol{q}_{n+1}(\boldsymbol{W}_h)+\boldsymbol{q}_n(\boldsymbol{W}_h)}{\boldsymbol{W}_h^{n+1}+\boldsymbol{W}_h^n}
	\nonumber\\=&\Delta t\sum_{n=0}^{m-1}\Lnorm{\boldsymbol{f}(t_{n+1})+\boldsymbol{f}(t_n)}{\boldsymbol{W}_h^{n+1}+\boldsymbol{W}_h^n}+\Delta t\sum_{n=0}^{m-1}\Gnorm{{\boldsymbol{g}_N(t_{n+1})}+{\boldsymbol{g}_N(t_n)}}{\boldsymbol{W}_h^{n+1}+\boldsymbol{W}_h^n}.\label{pf:frac:cg:T:sta:eq1}
	\end{align}	
	Expanding $\boldsymbol{q}_n$ allows us to rewrite \eqref{pf:frac:cg:T:sta:eq1} as
	\begin{align}
	2\rho&\llnorm{\boldsymbol{W}_h^{m}}^2+\frac{\Delta t^{2-\alpha}}{\Gamma(3-\alpha)}\sum_{n=0}^{m-1}\enorm{\boldsymbol{W}_h^{n+1}+\boldsymbol{W}_h^n}^2
	\nonumber\\
	=&2\rho\llnorm{\boldsymbol{W}_h^{0}}^2+\Delta t\sum_{n=0}^{m-1}\Lnorm{\boldsymbol{f}(t_{n+1})+\boldsymbol{f}(t_n)}{\boldsymbol{W}_h^{n+1}+\boldsymbol{W}_h^n}\nonumber\\&+\Delta t\sum_{n=0}^{m-1}\Gnorm{{\boldsymbol{g}_N(t_{n+1})}+{\boldsymbol{g}_N(t_n)}}{\boldsymbol{W}_h^{n+1}+\boldsymbol{W}_h^n}\nonumber\\&-\frac{\Delta t^{2-\alpha}}{\Gamma(3-\alpha)}\sum_{n=0}^{m-1}\Enorm{\sum_{i=0}^{n}B_{n+1,i}\boldsymbol{W}_h^i+\sum_{i=0}^{n-1}B_{n,i}\boldsymbol{W}_h^i}{\boldsymbol{W}_h^{n+1}+\boldsymbol{W}_h^n}.\label{pf:frac:cg:T:sta:eq2}
	\end{align}
	We shall find the bounds of the right hand side of \eqref{pf:frac:cg:T:sta:eq2}.
	\begin{itemize}[$\bullet$]
		\item $\llnorm{\boldsymbol{W}_h^{0}}^2$\\
		Since \eqref{eq:frac:T:fully:eq2} holds, we have
		\[\enorm{\boldsymbol W^0_h}^2\leq\enorm{\boldsymbol W^0_h}\enorm{\boldsymbol w_0}\]by Cauchy-Schwarz inequality and so \[\llnorm{\boldsymbol W^0_h}\leq\hnorm{\boldsymbol W^0_h}\leq C\enorm{\boldsymbol W^0_h}\leq C\enorm{\boldsymbol w_0}\]for some positive $C$ by norm equivalence between $H^1$ norm and energy norm.
		\item $ \Delta t\sum\limits_{n=0}^{m-1}\Lnorm{\boldsymbol{f}(t_{n+1})+\boldsymbol{f}(t_n)}{\boldsymbol{W}_h^{n+1}+\boldsymbol{W}_h^n} $\\
		Use of Cauchy-Schwarz and Young's inequalities gives
		\begin{align*}
		\Delta t&\sum\limits_{n=0}^{m-1}\Lnorm{\boldsymbol{f}(t_{n+1})+\boldsymbol{f}(t_n)}{\boldsymbol{W}_h^{n+1}+\boldsymbol{W}_h^n}\\
		\leq&\Delta t\sum\limits_{n=0}^{N}{2\epsilon_a}\llnorm{\boldsymbol{f}(t_n)}^2+\frac{2(T+\Delta t)}{\epsilon_a}\max_{0\leq n\leq N}\llnorm{\boldsymbol{W}^n_h}^2
		\end{align*}for any positive $\epsilon_a$.
		\item $ \Delta t\sum\limits_{n=0}^{m-1}\Gnorm{{\boldsymbol{g}_N(t_{n+1})}+{\boldsymbol{g}_N(t_n)}}{\boldsymbol{W}_h^{n+1}+\boldsymbol{W}_h^n} $\\
		While using the same approach as the above, we can also derive\begin{align*}
		\Delta t&\sum\limits_{n=0}^{m-1}\Gnorm{{\boldsymbol{g}_N(t_{n+1})}+{\boldsymbol{g}_N(t_n)}}{\boldsymbol{W}_h^{n+1}+\boldsymbol{W}_h^n}\\
		\leq&\Delta t\sum\limits_{n=0}^{N}{2\epsilon_b}\gnorm{{\boldsymbol{g}_N(t_n)}}^2+Ch^{-1}\frac{2(T+\Delta t)}{\epsilon_b}\max_{0\leq n\leq N}\llnorm{\boldsymbol{W}_h^n}^2
		\end{align*}by \eqref{eq:trace:inv}, for any positive $\epsilon_b$.
	\end{itemize}
	From the above bounds, \eqref{pf:frac:cg:T:sta:eq2} can be written as
	\begin{align}
	2\rho&\llnorm{\boldsymbol{W}_h^{m}}^2+\frac{\Delta t^{2-\alpha}}{\Gamma(3-\alpha)}\sum_{n=0}^{m-1}\enorm{\boldsymbol{W}_h^{n+1}+\boldsymbol{W}_h^n}^2
	\nonumber\\
	\leq&\rho C\enorm{\boldsymbol{w}_{0}}^2+\Delta t\sum\limits_{n=0}^{N}{2\epsilon_a}\llnorm{\boldsymbol{f}(t_n)}^2+\frac{2(T+\Delta t)}{\epsilon_a}\max_{0\leq n\leq N}\llnorm{\boldsymbol{W}^n_h}^2\nonumber\\&+\Delta t\sum\limits_{n=0}^{N}{2\epsilon_b}\gnorm{{\boldsymbol{g}_N(t_n)}}^2+\frac{Ch^{-1}2(T+\Delta t)}{\epsilon_b}\max_{0\leq n\leq N}\llnorm{\boldsymbol{W}_h^n}^2
	\nonumber\\&-\frac{\Delta t^{2-\alpha}}{\Gamma(3-\alpha)}\sum_{n=0}^{m-1}\Enorm{\sum_{i=0}^{n}B_{n+1,i}\boldsymbol{W}_h^i+\sum_{i=0}^{n-1}B_{n,i}\boldsymbol{W}_h^i}{\boldsymbol{W}_h^{n+1}+\boldsymbol{W}_h^n}\nonumber
	\\
	:=&\mathcal R-\frac{\Delta t^{2-\alpha}}{\Gamma(3-\alpha)}\sum_{n=0}^{m-1}\Enorm{\sum_{i=0}^{n}B_{n+1,i}\boldsymbol{W}_h^i+\sum_{i=0}^{n-1}B_{n,i}\boldsymbol{W}_h^i}{\boldsymbol{W}_h^{n+1}+\boldsymbol{W}_h^n}.\label{pf:frac:cg:T:sta:eq3}
	\end{align}Now, it remains to show the boundedness of \eqref{pf:frac:cg:T:sta:eq3}. Note that $\mathcal R$ in \eqref{pf:frac:cg:T:sta:eq3} is independent of $m$. Hereafter, we would like to use mathematical induction to derive the upper bound of the last term. Our claim to be shown by induction is \begin{align}
	2\rho&\llnorm{\boldsymbol{W}_h^{m}}^2+\frac{\Delta t^{2-\alpha}}{2\Gamma(3-\alpha)}\sum_{n=0}^{m-1}\enorm{\boldsymbol{W}_h^{n+1}+\boldsymbol{W}_h^n}^2
	\leq C\big(\mathcal R+{\Delta t^{2-\alpha}}\enorm{\boldsymbol{W}_h^{0}}^2\big),\label{pf:frac:cg:T:sta:eq5}
	\end{align}for some positive $C$, $\forall m$. For $m=1$ in the last term of \eqref{pf:frac:cg:T:sta:eq3}, we have
	\begin{align*}|\Enorm{B_{1,0}\boldsymbol{W}_h^0}{\boldsymbol{W}_h^{1}+\boldsymbol{W}_h^0}|\leq\frac{B_{1,0}^2\epsilon}{2}\enorm{\boldsymbol{W}_h^0}^2+\frac{1}{2\epsilon}\enorm{\boldsymbol{W}_h^{1}+\boldsymbol{W}_h^0}^2
	\end{align*}by Cauchy-Schwarz and Young's inequalities with any positive $\epsilon$. Hence, taking $\epsilon=1$ allows us to have
	\begin{align}
	2\rho&\llnorm{\boldsymbol{W}_h^{1}}^2+\frac{\Delta t^{2-\alpha}}{2\Gamma(3-\alpha)}\enorm{\boldsymbol{W}_h^{1}+\boldsymbol{W}_h^0}^2
	\leq\mathcal R+\frac{B_{1,0}^2}{2}\frac{\Delta t^{2-\alpha}}{\Gamma(3-\alpha)}\enorm{\boldsymbol{W}_h^{0}}^2.\label{pf:frac:cg:induction:m1}
	\end{align}
	When $m=2$, \eqref{pf:frac:cg:T:sta:eq3} gives
	\begin{align*}
	2\rho&\llnorm{\boldsymbol{W}_h^{2}}^2+\frac{\Delta t^{2-\alpha}}{\Gamma(3-\alpha)}(\enorm{\boldsymbol{W}_h^{2}+\boldsymbol{W}_h^1}^2+\enorm{\boldsymbol{W}_h^{1}+\boldsymbol{W}_h^0}^2)\\
	\leq&\mathcal R-\frac{\Delta t^{2-\alpha}}{\Gamma(3-\alpha)}\Enorm{B_{2,1}\boldsymbol{W}_h^1+B_{2,0}\boldsymbol{W}_h^0+B_{1,0}\boldsymbol{W}_h^0}{\boldsymbol{W}_h^{2}+\boldsymbol{W}_h^1}\\&-\frac{\Delta t^{2-\alpha}}{\Gamma(3-\alpha)}\Enorm{B_{1,0}\boldsymbol{W}_h^0}{\boldsymbol{W}_h^{1}+\boldsymbol{W}_h^0}.
	\end{align*}Using Cauchy-Schwarz and Young's inequalities, we can have
	\begin{align*}
	\big|&\Enorm{B_{2,1}\boldsymbol{W}_h^1+B_{2,0}\boldsymbol{W}_h^0+B_{1,0}\boldsymbol{W}_h^0}{\boldsymbol{W}_h^{2}+\boldsymbol{W}_h^1}\big|\\
	&\leq\frac{(\max(B_{2,1},B_{2,0}+B_{1,0}))^2}{2\epsilon}\enorm{\boldsymbol{W}_h^{1}+\boldsymbol{W}_h^{0}}^2+\frac{\epsilon}{2}\enorm{\boldsymbol{W}_h^{2}+\boldsymbol{W}_h^{1}}^2,
	\end{align*}and\begin{align*}
	\big|\Enorm{B_{1,0}\boldsymbol{W}_h^0}{\boldsymbol{W}_h^{1}+\boldsymbol{W}_h^0}
\big|	\leq\frac{B_{1,0}^2}{2\epsilon}\enorm{\boldsymbol{W}_h^{0}}^2+\frac{\epsilon}{2}\enorm{\boldsymbol{W}_h^{1}+\boldsymbol{W}_h^{0}}^2
	\end{align*}for any positive $\epsilon$. Hence coupling with \eqref{pf:frac:cg:induction:m1} which provides the bound for $\enorm{\boldsymbol{W}_h^{1}+\boldsymbol{W}_h^{0}}^2$, and choosing $\epsilon=1$, we can write \eqref{pf:frac:cg:T:sta:eq3} for $m=2$ as
	\begin{align*}
	2\rho&\llnorm{\boldsymbol{W}_h^{2}}^2+\frac{\Delta t^{2-\alpha}}{2\Gamma(3-\alpha)}\sum_{n=0}^{1}\enorm{\boldsymbol{W}_h^{n+1}+\boldsymbol{W}_h^n}^2\leq C\big(\mathcal R+\Delta t^{2-\alpha}\enorm{\boldsymbol{W}_h^{0}}^2\big),
	\end{align*}for some positive $C$. Let us assume that \eqref{pf:frac:cg:T:sta:eq5} holds for $m=j<N$ so that
	\begin{align*}
	2\rho\llnorm{\boldsymbol{W}_h^{j}}^2+\frac{\Delta t^{2-\alpha}}{2\Gamma(3-\alpha)}\sum_{n=0}^{j-1}\enorm{\boldsymbol{W}_h^{n+1}+\boldsymbol{W}_h^n}^2
	\leq&C\big(\mathcal R+\Delta t^{2-\alpha}\enorm{\boldsymbol{W}_h^{0}}^2\big),
	\end{align*}and then for $m=j+1$, we have, from \eqref{pf:frac:cg:T:sta:eq3},
	\begin{align*}\sum_{n=0}^{j}&\Enorm{\sum_{i=0}^{n}B_{n+1,i}\boldsymbol{W}_h^i+\sum_{i=0}^{n-1}B_{n,i}\boldsymbol{W}_h^i}{\boldsymbol{W}_h^{n+1}+\boldsymbol{W}_h^n}\\
	\leq&\sum_{n=0}^{j}\sum_{i=1}^{n-1}\left(\frac{G^2\epsilon}{2}\enorm{\boldsymbol{W}_h^{i+1}+\boldsymbol{W}_h^i}^2+\frac{1}{2\epsilon}\enorm{\boldsymbol{W}_h^{n+1}+\boldsymbol{W}_h^n}^2\right)\\&+\sum_{n=0}^{j}\left(\frac{(3G)^2\tilde\epsilon}{2}\enorm{\boldsymbol{W}_h^0}^2+\frac{1}{2\tilde\epsilon}\enorm{\boldsymbol{W}_h^{n+1}+\boldsymbol{W}_h^n}^2\right)\\&+\sum_{n=0}^{j}\left(\frac{G^2\check\epsilon}{2}\enorm{\boldsymbol{W}_h^1+\boldsymbol{W}_h^0}^2+\frac{1}{2\check\epsilon}\enorm{\boldsymbol{W}_h^{n+1}+\boldsymbol{W}_h^n}^2\right)
	\end{align*}where $0<G=\max\limits_{0\leq i\leq n\leq N}B_{n,i}<2$, for any positive $\epsilon,\tilde\epsilon,\text{ and }\check\epsilon$. Since $\sum\limits_{i=1}^{n-1}\DGnorm{\boldsymbol{W}_h^{i+1}+\boldsymbol{W}_h^i}^2$ is bounded for \mbox{$0\leq n\leq j$} by the induction assumption, we can obtain the boundedness of $\sum\limits_{n=0}^{j}\sum\limits_{i=1}^{n-1}\enorm{\boldsymbol{W}_h^{i+1}+\boldsymbol{W}_h^i}^2$. Consequently, setting $\epsilon=\tilde\epsilon=\check\epsilon=1/3$ yields
	\begin{align*}
	2\rho&\llnorm{\boldsymbol{W}_h^{j+1}}^2+\frac{\Delta t^{2-\alpha}}{2\Gamma(3-\alpha)}\sum_{n=0}^{j}\enorm{\boldsymbol{W}_h^{n+1}+\boldsymbol{W}_h^n}^2
	\leq C\bigg(\mathcal R+{\Delta t^{2-\alpha}}\enorm{\boldsymbol{W}_h^{0}}^2\bigg).
	\end{align*}Thus we can complete the induction and hence \eqref{pf:frac:cg:T:sta:eq5} holds. Turning to our main goal, when we consider maximum in \eqref{pf:frac:cg:T:sta:eq5} with the argument
	\[a_n+b_n\leq C\Rightarrow \max_na_n+\max_nb_n\leq 2C, \text{ for any positive }a_n,b_n,\forall n,\]then \eqref{pf:frac:cg:T:sta:eq5} can be written as
	\begin{align*}
	2\rho&\max_{0\leq n\leq N}\llnorm{\boldsymbol{W}_h^{n}}^2+\frac{\Delta t^{2-\alpha}}{2\Gamma(3-\alpha)}\sum_{n=0}^{N-1}\enorm{\boldsymbol{W}_h^{n+1}+\boldsymbol{W}_h^n}^2
	\\
	\leq&2C\bigg(\enorm{\boldsymbol{w}_{0}}^2+\Delta t\sum\limits_{n=0}^{N}{2\epsilon_a}\llnorm{\boldsymbol{f}(t_n)}^2+\frac{2(T+\Delta t)}{\epsilon_a}\max_{0\leq n\leq N}\llnorm{\boldsymbol{W}^n_h}^2\nonumber\\&+\Delta t\sum\limits_{n=0}^{N}{2\epsilon_b}\gnorm{{\boldsymbol{g}_N(t_n)}}^2+\frac{h^{-1}2(T+\Delta t)}{\epsilon_b}\max_{0\leq n\leq N}\llnorm{\boldsymbol{W}_h^n}^2\bigg),
	\end{align*}since $\enorm{\boldsymbol{W}_h^{0}}\leq\enorm{\boldsymbol{w}_{0}}$.
	Therefore, choosing $\epsilon_a=8C(T+\Delta t)/\rho$ and $\epsilon_b=8Ch^{-1}(T+\Delta t)/\rho$ leads us to have
	\begin{align*}
	\rho&\max_{0\leq n\leq N}\llnorm{\boldsymbol{W}_h^{n}}^2+\frac{\Delta t^{2-\alpha}}{2\Gamma(3-\alpha)}\sum_{n=0}^{N-1}\enorm{\boldsymbol{W}_h^{n+1}+\boldsymbol{W}_h^n}^2
	\nonumber\\
	\leq&CT\bigg(\enorm{\boldsymbol{w}_{0}}^2+\Delta t\sum\limits_{n=0}^{N}\llnorm{\boldsymbol{f}(t_n)}^2+\Delta t\sum\limits_{n=0}^{N}h^{-1}\gnorm{{\boldsymbol{g}_N(t_n)}}^2\bigg).
	\end{align*}Furthermore, this bound also implies the existence and uniqueness of the discrete solution.
\end{proof}

\begin{rmk}
	In Theorem \ref{thm:frac:cg:fully:sta}, $h^{-1}$ term appears on the traction. However, in practice, it has nothing to do because of fixed $h$ for the finite element spaces. It comes from the fact that the boundary condition is imposed weakly \cite{DGV,DG} on account of inverse polynomial trace theorem.
\end{rmk}
\begin{rmk}
	If we use Gr\"onwall's inequality to show the boundedness rather than taking the maximum, we have an exponentially increasing bound in the final time $T$, e.g. see \cite{li2011developing,thomee1984galerkin}. That is, instead of $T$, we have $\exp(T)$ on the stability bound. 
\end{rmk}\begin{rmk}
In Theorem \ref{thm:frac:cg:semi:sta}, the stability bound has been proved by the positive definiteness of the kernel (fractional integral). On the other hand, our discrete kernel is no longer positive definite but weak positive definite. We refer \cite{mclean1993numerical,grenander1958toeplitz,lopez1990difference} for more details. Moreover, in order to use the positive definiteness for the stability analysis, zero traction or pure Dirichlet boundary condition should be further assumed. However, we employ inverse trace polynomial theorem for mixed boundary conditions in the proof of Theorem \ref{thm:frac:cg:fully:sta} with by means of induction rather than the use of positive definiteness. 
\end{rmk}
\subsection{Error Estimates}\label{subsec:error}
\noindent In terms of errors analysis in time, the Crank-Nicolson method requires at least $H^3$ smoothness of a solution with respect to time to get second order accuracy. However, due to the weak singularity, we may encounter restrictions on high regularity of solutions. Hence we will remark on the regularity of solutions.
\begin{rmk}(Regularity of solutions)		
	Let us recall the primal equation \eqref{eq:frac:primal:eq1}. We can rewrite it in convolution form so that
	\begin{align*}
	\rho\dot{\boldsymbol{w}}(t)=&\nabla\cdot\leftidx{_0}I^{1-\alpha}_t({\hat{\ushort{\boldsymbol{D}}}}\ushort{{\boldsymbol{\varepsilon}}}(\boldsymbol{w}(t)))+\boldsymbol{f}(t)\\
	=&\beta_{1-\alpha}*\mathcal{D}\boldsymbol{w}(t)+\boldsymbol{f}(t)
	\end{align*}where $\mathcal{D}=\nabla\cdot {\hat{\ushort{\boldsymbol{D}}}}\ushort{{\boldsymbol{\varepsilon}}}$ is a linear differential operator on the spatial domain and $\beta_{1-\alpha}=t^{-\alpha}/\Gamma(1-\alpha)$ for $\alpha\in(0,1)$. By Young's inequality for the convolution, we can observe that\begin{align*}
	\lVert\rho\dot{\boldsymbol{w}}\rVert_{L_q(0,T)}\leq&\lVert\beta_{1-\alpha}\rVert_{L_1(0,T)}\lVert\mathcal{D}\boldsymbol{w}\rVert_{L_q(0,T)}+\lVert\boldsymbol{f}\rVert_{L_q(0,T)}\text{ for }q\geq1.
	\end{align*}Since $\beta_{1-\alpha}$ is $L_1$ integrable, if $\mathcal{D}\boldsymbol{w}$ and $\boldsymbol{f}$ are $L_2$ integrable in time, so is $\dot{\boldsymbol{w}}$. Differentiating \eqref{eq:frac:primal:eq1} with respect to time gives
	\begin{align*}
	\rho\ddot{\boldsymbol{w}}(t)
	=&\beta_{1-\alpha}(t)\mathcal{D}\boldsymbol{w}(0)+\beta_{1-\alpha}*\mathcal{D}\dot{\boldsymbol{w}}(t)+\dot{\boldsymbol{f}}(t).	
	\end{align*}
	We assume that $\boldsymbol{w}(0)\in [H^2(\Omega)]^d$ then $\ddot{\boldsymbol{w}}$ is $L_1$ integrable with $L_1$ integrable $\dot{\boldsymbol{f}}$ and $\mathcal{D}\dot{\boldsymbol{w}}$ with respect to time. In this manner, we can observe $L_2$ integrable $\ddot{\boldsymbol{w}}$ if $\dot{\boldsymbol{f}}$ and $\mathcal{D}\dot{\boldsymbol{w}}$ are in $L_2$, and $\mathcal{D}\boldsymbol{w}(0)=\boldsymbol{0}$. Repeatedly, we can consider the third time derivative of $\boldsymbol{w}$. Then we have
	\begin{align}
	\rho{\boldsymbol{w}}^{(3)}(t)
	=&\dot\beta_{1-\alpha}(t)\mathcal{D}\boldsymbol{w}(0)+\beta_{1-\alpha}(t)\mathcal{D}\dot{\boldsymbol{w}}(0)+\beta_{1-\alpha}*\mathcal{D}\ddot{\boldsymbol{w}}(t)+\ddot{\boldsymbol{ f}}(t)	.
	\label{eq:reg:conv}\end{align}Note that $\dot\beta_{1-\alpha}(t)$ is non-integrable in $L_1$ and $L_2$ so it is not obviously seen that the third time derivative of $\boldsymbol{w}$ is integrable. In a second order finite difference scheme, the third derivative and its boundedness are required to take a full advantage of the second order schemes. For example, we can observe that if $\boldsymbol w$ is three-times differentiable,
	\begin{align}
	\frac{\dot {\boldsymbol w}(t_{n+1})+\dot{\boldsymbol  w}(t_{n})}{2}-\frac{\boldsymbol w(t_{n+1})-\boldsymbol w(t_{n})}{\Delta t}=\frac{1}{2\Delta t}\int^{t_{n+1}}_{t_n}\boldsymbol w^{(3)}(t)(t_{n+1}-t)(t-t_n)dt.\label{eq:fdm:est}
	\end{align}Moreover, the boundedness of $\boldsymbol w^{(3)}$ leads that \eqref{eq:fdm:est} is of order $\Delta t^2$.
	For example, when we suppose $\boldsymbol w^{(3)}\in{L_2(t_n,t_{n+1};L_2(\Omega))}$, the use of Cauchy-Schwarz inequality gives
	\begin{equation}
	\llnorm{\frac{\dot{\boldsymbol  w}(t_{n+1})+\dot {\boldsymbol w}(t_{n})}{2}-\frac{\boldsymbol w(t_{n+1})-\boldsymbol w(t_{n})}{\Delta t}}^2\leq\frac{\Delta t^3}{4}\lVert \boldsymbol w^{(3)}\rVert_{L_2(t_n,t_{n+1};L_2(\Omega))}^2.	\label{eq:frac:CN}\end{equation}
	By substitution of \eqref{eq:reg:conv} into \eqref{eq:fdm:est}, we can also observe that
	\begin{align}
	&\frac{\dot {\boldsymbol w}(t_{n+1})+\dot{\boldsymbol  w}(t_{n})}{2}-\frac{\boldsymbol w(t_{n+1})-\boldsymbol w(t_{n})}{\Delta t}\nonumber\\&=\frac{1}{2\rho\Delta t}\int^{t_{n+1}}_{t_n}\left(\dot\beta_{1-\alpha}(t)\mathcal{D}\boldsymbol{w}(0)+\beta_{1-\alpha}(t)\mathcal{D}\dot{\boldsymbol{w}}(0)+\beta_{1-\alpha}*\mathcal{D}\ddot{\boldsymbol{w}}(t)+\ddot{\boldsymbol{ f}}(t)\right)(t_{n+1}-t)(t-t_n)dt.\label{eq:fdm:est2}
	\end{align}
	Note that we assume $\boldsymbol{w}(0)\in[H^2(\Omega)]^d$ and $\dot{\boldsymbol{w}}(0)\in[H^2(\Omega)]^d$.
	Thus, if $\ddot{\boldsymbol{ f}}\in L_2(t_n,t_{n+1};[L_2(\Omega)]^d)$ and $\ddot{\boldsymbol{w}}\in L_2(t_n,t_{n+1};[H^2(\Omega)]^d)$, we have $ \boldsymbol w^{(3)}\in{L_2(t_n,t_{n+1};[L_2(\Omega)]^d)}$ for $n\geq1$.
	However, the singularity appears for $n=0$. So, we need to introduce the following lemma.	
\end{rmk}
	\begin{lemma}\label{lemma:regularity}
		Suppose  $\boldsymbol{f}\in W^3_1(0,T;[L_2(\Omega)]^d)\cap H^2(0,T;[H^2(\Omega)]^d)$, ${\boldsymbol{w}}(0)\in [H^2(\Omega)]^d$ and $\dot{\boldsymbol{w}}(0)\in [H^2(\Omega)]^d$. If ${\boldsymbol{w}}(0)={\boldsymbol{0}}$, we have
\begin{equation}
\llnorm{\frac{\boldsymbol{w}(t_1)-\boldsymbol{w}(t_0)}{\Delta t}-\frac{\dot {\boldsymbol{w}}(t_1)+\dot{\boldsymbol{w}} (t_0)}{2}}\leq C  \Delta t^{2-\alpha},\label{eq:frac:CN:t0}
\end{equation}for some positive constant $C$ independent of $\Delta t$. Furthermore, we can also obtain
\begin{equation}
\llnorm{\frac{\boldsymbol{w}(t_1)-\boldsymbol{w}(t_0)}{\Delta t}-\frac{\dot {\boldsymbol{w}}(t_1)+\dot{\boldsymbol{w}} (t_0)}{2}}\leq C  \Delta t^{2},\label{eq:frac:CN:t0:opt}
\end{equation}when ${\boldsymbol{w}}(0)=\dot{\boldsymbol{w}}(0)={\boldsymbol{0}}$.
	\end{lemma}
	\begin{proof}
		Let us recall \eqref{eq:fdm:est2}
		\begin{align*}
		&\frac{\dot {\boldsymbol w}(t_{n+1})+\dot{\boldsymbol  w}(t_{n})}{2}-\frac{\boldsymbol w(t_{n+1})-\boldsymbol w(t_{n})}{\Delta t}\\&=\frac{1}{2\rho\Delta t}\int^{t_{n+1}}_{t_n}\left(\dot\beta_{1-\alpha}(t)\mathcal{D}\boldsymbol{w}(0)+\beta_{1-\alpha}(t)\mathcal{D}\dot{\boldsymbol{w}}(0)+\beta_{1-\alpha}*\mathcal{D}\ddot{\boldsymbol{w}}(t)+\ddot{\boldsymbol{f}}(t)\right)(t_{n+1}-t)(t-t_n)dt.
		\end{align*}
		We can expand it by
		\begin{align*}
		\frac{\dot {\boldsymbol w}(t_{n+1})+\dot{\boldsymbol  w}(t_{n})}{2}-\frac{\boldsymbol w(t_{n+1})-\boldsymbol w(t_{n})}{\Delta t}=&\frac{1}{2\rho\Delta t}\bigg(\int^{t_{n+1}}_{t_n}\frac{-\alpha t^{-\alpha-1}}{\Gamma(1-\alpha)}\mathcal{D}\boldsymbol{w}(0)(t_{n+1}-t)(t-t_n)dt\\
		&+\int^{t_{n+1}}_{t_n}\frac{ t^{-\alpha}}{\Gamma(1-\alpha)}\mathcal{D}\dot{\boldsymbol{w}}(0)(t_{n+1}-t)(t-t_n)dt\\
		&+\int^{t_{n+1}}_{t_n}\left(\beta_{1-\alpha}*\mathcal{D}\ddot{\boldsymbol{w}}(t)+\ddot{\boldsymbol{ f}}(t)\right)(t_{n+1}-t)(t-t_n)dt\bigg).
		\end{align*}
		Consider the first and second terms of the right hand side for $n=0$. Then we have
		\begin{align*}
		\int^{\Delta t}_{0}\frac{-\alpha t^{-\alpha-1}}{\Gamma(1-\alpha)}\mathcal{D}\boldsymbol{w}(0)(\Delta t-t)tdt=\frac{-\alpha \Delta t^{2-\alpha}}{\Gamma(3-\alpha)}\mathcal{D}\boldsymbol{w}(0),
		\end{align*}
		and
		\begin{align*}
		\int^{\Delta t}_{0}\frac{ t^{-\alpha}}{\Gamma(1-\alpha)}\mathcal{D}\dot{\boldsymbol{w}}(0)(\Delta t-t)tdt=\frac{(1-\alpha) \Delta t^{3-\alpha}}{\Gamma(4-\alpha)}\mathcal{D}\dot{\boldsymbol{w}}(0).
		\end{align*}
		Thus, Cauchy-Schwarz inequality and Young's inequality for convolution lead us to have
		\begin{align*}
		\llnorm{
		\frac{\dot {\boldsymbol w}(t_{1})+\dot{\boldsymbol  w}(t_{0})}{2}-\frac{\boldsymbol w(t_{1})-\boldsymbol w(t_{0})}{\Delta t}}\leq&\frac{\alpha\Delta t^{1-\alpha}}{2\rho\Gamma(3-\alpha)}\llnorm{\mathcal{D}{\boldsymbol{w}}(0)}+\frac{(1-\alpha)\Delta t^{2-\alpha}}{2\rho\Gamma(4-\alpha)}\llnorm{\mathcal{D}\dot{\boldsymbol{w}}(0)}\\&+C{\Delta t^2},
		\end{align*}where $C$ depends on $f$ but is independent of $\Delta t$, e.g. see \cite{mclean1993numerical} for more details.
		We can conclude that if $\boldsymbol{w}(0)=\boldsymbol{0}$, \eqref{eq:frac:CN:t0} holds. Moreover, when we additionally assume $\dot{\boldsymbol{w}}(0)=\boldsymbol{0}$, we obtain
		\begin{equation*}
		\llnorm{\frac{\boldsymbol{w}(t_1)-\boldsymbol{w}(t_0)}{\Delta t}-\frac{\dot {\boldsymbol{w}}(t_1)+\dot{\boldsymbol{w}} (t_0)}{2}}\leq C  \Delta t^{2}.
		\end{equation*}
	\end{proof}
	
	\begin{rmk}
		We refer to \cite{mclean1993numerical} for the assumption $\boldsymbol{f}\in W^3_1(0,T;[L_2(\Omega)]^d)\cap H^2(0,T;[H^2(\Omega)]^d)$. In addition, once $\boldsymbol w\in\mathrm{ker}(\mathcal{D})$ 
		where $\mathrm{ker}(\mathcal{D})$ is a kernel set of the differential operator $\mathcal{D}$, the strong form becomes a simple first order ODE problem so that the singularity will also disappear.		
	\end{rmk} In order to consider spatial error estimates, we want to introduce the following elliptic error estimates. We define an elliptic projection $\boldsymbol{R}:\boldsymbol{V}\mapsto\boldsymbol{V}^h$ by
\[\Enorm{\boldsymbol{R}\boldsymbol{w}}{\boldsymbol{v}}=\Enorm{\boldsymbol{w}}{\boldsymbol{v}},\text{for ${\boldsymbol{w}}\in\boldsymbol{V}$ and any ${\boldsymbol{v}}\in\boldsymbol{V}^h$},\]then we have \textit{Galerkin orthogonality} such that for ${\boldsymbol{w}}\in\boldsymbol{V}$,
\[\Enorm{\boldsymbol{R}\boldsymbol{w}-\boldsymbol{w}}{\boldsymbol{v}}=0,\qquad\text{for any ${\boldsymbol{v}}\in\boldsymbol{V}^h$}.\]
According to \cite{MTE,wheeler}, we can obtain elliptic error estimates such that
\begin{align}
\enorm{\boldsymbol w-\boldsymbol R\boldsymbol w}\leq C|\boldsymbol w|_{H^{r}(\Omega)}h^{r-1},\label{eq:elliptic:energy}
\end{align}where $\boldsymbol V^h\subset\boldsymbol V$ is a subspace of polynomials of degree $k$, $\boldsymbol w\in [H^{s}(\Omega)]^d$, and \mbox{$r=\min(k+1,s)$.} Moreover, the use of elliptic regularity estimates \cite{MTE,dauge1988elliptic,grisvard2011elliptic} in a standard duality argument enables us to get\begin{align}
\llnorm{\boldsymbol w-\boldsymbol R\boldsymbol w}\leq C|\boldsymbol w|_{H^{r}(\Omega)}h^{r}\label{eq:elliptic:l2}.
\end{align}

Next, we state and prove \textit{a priori} error estimates by recalling elliptic approximations \eqref{eq:elliptic:energy} and \eqref{eq:elliptic:l2}. Hence we use the elliptic projection operator $\boldsymbol{R}$ and define\[\boldsymbol \theta(t):=\boldsymbol w(t)-\boldsymbol R\boldsymbol w(t),\ \text{for }t\in[0,T],\qquad
\boldsymbol \chi^n:=\boldsymbol W_h^n-\boldsymbol R\boldsymbol w(t_n)\textrm{
	for $n=0,\ldots,N$}.\] \begin{lemma}\label{lemma:frac:error:T}
	Suppose $\boldsymbol{f}$ and $\boldsymbol{w}_0$ are given to hold \eqref{eq:frac:CN:t0}, \[\boldsymbol{w}\in C^2(0,T;[H^{s}(\Omega)]^d)\cap W^1_\infty(0,T;\boldsymbol{V})\text{ for $s\geq2$},\] and $\left(\boldsymbol{W}^n_h\right)_{n=0}^N$ satisfies the fully discrete formulae \eqref{eq:frac:T:fully:eq1} and \eqref{eq:frac:T:fully:eq2}. Then we have
	\begin{align*}
	&\max_{0\leq n\leq N}\llnorm{\boldsymbol{\chi}^n}+\left(\Delta t^{2-\alpha}\sum_{n=0}^{N-1}\enorm{\boldsymbol{\chi}^{n+1}+\boldsymbol{\chi}^n}^2\right)^{1/2}\leq CT^{2-\alpha}(h^{r}+\Delta t^{2-\alpha}),
	\end{align*}where positive $C$ is independent of $h$ and $\Delta t$, and
	$r=\min(k+1,s)$. Moreover, if $\boldsymbol{w}(0)=\dot{\boldsymbol{w}}(0)=\boldsymbol{0}$ or $\boldsymbol{w}\in\ker(\mathcal{D})$, then we have
	\begin{align*}
	&\max_{0\leq n\leq N}\llnorm{\boldsymbol{\chi}^n}+\left(\Delta t^{2-\alpha}\sum_{n=0}^{N-1}\enorm{\boldsymbol{\chi}^{n+1}+\boldsymbol{\chi}^n}^2\right)^{1/2}= CT^{2-\alpha}(h^{r}+\Delta t^{2}).
	\end{align*}
\end{lemma}
\begin{proof}
	For $m\in\{1,\ldots,N\}$, subtracting the average of \eqref{eq:frac:T:primal:eq1} over $t=t_{n+1}$ and $t=t_n$ from \eqref{eq:frac:T:fully:eq1} where $0\leq n\leq m-1$ gives
	\begin{align*}
	\rho\Lnorm{\frac{\boldsymbol{W}_h^{n+1}-\boldsymbol{W}_h^{n}}{\Delta t}-\frac{\dot{ \boldsymbol w}^{n+1}+\dot{\boldsymbol w}^{n}}{2}}{\boldsymbol{v}}+\Enorm{\frac{\boldsymbol{q}_{n+1}(\boldsymbol{W}_h)+\boldsymbol{q}_{n}(\boldsymbol{W}_h)}{2}-\frac{\leftidx_{0}I^{1-\alpha}_{t_{n+1}}\boldsymbol{{w}}+\leftidx_{0}I^{1-\alpha}_{t_{n}}\boldsymbol{{w}}}{2}}{\boldsymbol{v}}=0
	\end{align*}for any $\boldsymbol{v}\in\boldsymbol{V}^h$.
	By definitions of $\boldsymbol\theta$ and $\boldsymbol\chi$, we can rewrite this as
	\begin{align}
	\frac{\rho}{\Delta t}&\Lnorm{\boldsymbol{\chi}^{n+1}-\boldsymbol{\chi}^{n}}{\boldsymbol{v}}+\frac{1}{2}\Enorm{\boldsymbol{q}_{n+1}(\boldsymbol{\chi})+\boldsymbol{q}_{n}(\boldsymbol{\chi})}{\boldsymbol{v}}\nonumber\\=&\frac{\rho}{\Delta t}\Lnorm{\boldsymbol{\theta}^{n+1}-\boldsymbol{\theta}^{n}}{\boldsymbol{v}}+\frac{1}{2}\Enorm{\boldsymbol{q}_{n+1}(\boldsymbol{\theta})+\boldsymbol{q}_{n}(\boldsymbol{\theta})}{\boldsymbol{v}}+\frac{1}{2}\Enorm{\boldsymbol{e}^{n+1}+\boldsymbol{e}^n}{\boldsymbol{v}}+\rho\Lnorm{\boldsymbol {\mathcal{E}}^n}{\boldsymbol{v}},\label{pf:frac:lemma:error:T:eq1}
	\end{align} where $\boldsymbol{e}^n:=\boldsymbol{q}_{n}(\boldsymbol{w})-\leftidx_{0}I^{1-\alpha}_{t}\boldsymbol{w}(t_n)$ and $\boldsymbol {\mathcal{E}}(t):=\frac{\dot {\boldsymbol w}(t+\Delta t)+\dot {\boldsymbol w}(t)}{2}-\frac{\boldsymbol w(t+\Delta t)-\boldsymbol w(t)}{\Delta t}$ for $t\in[0,T-\Delta t]$. Galerkin orthogonality reduces \eqref{pf:frac:lemma:error:T:eq1} to \begin{align}
	\frac{\rho}{\Delta t}\Lnorm{\boldsymbol{\chi}^{n+1}-\boldsymbol{\chi}^{n}}{\boldsymbol{v}}+\frac{1}{2}\Enorm{\boldsymbol{q}_{n+1}(\boldsymbol{\chi})+\boldsymbol{q}_{n}(\boldsymbol{\chi})}{\boldsymbol{v}}=\frac{\rho}{\Delta t}\Lnorm{\boldsymbol{\theta}^{n+1}-\boldsymbol{\theta}^{n}}{\boldsymbol{v}}+\frac{1}{2}\Enorm{\boldsymbol{e}^{n+1}+\boldsymbol{e}^n}{\boldsymbol{v}}+\rho\Lnorm{\boldsymbol {\mathcal{E}}^n}{\boldsymbol{v}}.\label{pf:frac:lemma:error:T:eq2}
	\end{align}
	Once we put  $\boldsymbol{v}=2\Delta t(\boldsymbol{\chi}^{n+1}+\boldsymbol{\chi}^{n})$ in \eqref{pf:frac:lemma:error:T:eq2}, summing from $n=0$ to $n=m-1$ produces
	\begin{align}
	2\rho&\llnorm{\boldsymbol{\chi}^m}^2+\frac{\Delta t^{2-\alpha}}{\Gamma(3-\alpha)}\sum_{n=0}^{m-1}\enorm{\boldsymbol{\chi}^{n+1}+\boldsymbol{\chi}^n}^2
	\nonumber\\
	=&2\rho\llnorm{\boldsymbol{\chi}^0}^2+2{\rho}\sum_{n=0}^{m-1}\Lnorm{\boldsymbol{\theta}^{n+1}-\boldsymbol{\theta}^{n}}{{\boldsymbol{\chi}^{n+1}+\boldsymbol{\chi}^n}}+\Delta t\sum_{n=0}^{m-1}\Enorm{\boldsymbol{e}^{n+1}+\boldsymbol{e}^n}{{\boldsymbol{\chi}^{n+1}+\boldsymbol{\chi}^n}}\nonumber\\&+2\rho\Delta t\sum_{n=0}^{m-1}\Lnorm{\boldsymbol {\mathcal{E}}^n}{{\boldsymbol{\chi}^{n+1}+\boldsymbol{\chi}^n}}-\frac{\Delta t^{2-\alpha}}{\Gamma(3-\alpha)}\sum_{n=0}^{m-1}\Enorm{\sum_{i=0}^{n}B_{n+1,i}\boldsymbol{\chi}^i+\sum_{i=0}^{n-1}B_{n,i}\boldsymbol{\chi}^i}{\boldsymbol{\chi}^{n+1}+\boldsymbol{\chi}^n}.\label{pf:frac:lemma:error:T:eq3}
	\end{align}	For the sake of error estimation, we shall show the bounds of \eqref{pf:frac:lemma:error:T:eq3} as following.
	\begin{itemize}[$\bullet$]
		\item $\llnorm{\boldsymbol{\chi}^0}^2$\\
		\eqref{eq:frac:T:fully:eq2} and Galerkin orthogonality lead us to have \begin{align*}
		\Enorm{\boldsymbol{\chi}^0}{\boldsymbol{v}}
		=&\Enorm{\boldsymbol{W}_h^0+(\boldsymbol{w}_0-\boldsymbol{w}_0)-\boldsymbol{R}\boldsymbol{w}_0}{\boldsymbol{v}}=
		\Enorm{\boldsymbol{W}_h^0-\boldsymbol{w}_0}{\boldsymbol{v}}+
		\Enorm{\boldsymbol{w}_0-\boldsymbol{R}\boldsymbol{w}_0}{\boldsymbol{v}}
		=0\end{align*}for any $\boldsymbol{v}\in\boldsymbol{V}^h$ and hence
		$\llnorm{\boldsymbol{\chi}^0}^2=0$.
		\item $\sum\limits_{n=0}^{m-1}\Lnorm{\boldsymbol{\theta}^{n+1}-\boldsymbol{\theta}^{n}}{{\boldsymbol{\chi}^{n+1}+\boldsymbol{\chi}^n}}$\\
		Since $\boldsymbol{w}$ belongs to $H^1$ in time, we can write
		\begin{align*}
		\sum\limits_{n=0}^{m-1}&\Lnorm{\boldsymbol{\theta}^{n+1}-\boldsymbol{\theta}^{n}}{{\boldsymbol{\chi}^{n+1}+\boldsymbol{\chi}^n}}
		=\sum\limits_{n=0}^{m-1}\int_{t_{n}}^{t_{n+1}}\Lnorm{\dot{\boldsymbol{\theta}}(t')}{{\boldsymbol{\chi}^{n+1}+\boldsymbol{\chi}^n}}dt'\\
		\leq&\frac{\epsilon_a}{2}\int_{0}^{t_{m}}\llnorm{\dot{\boldsymbol{\theta}}(t')}^2dt'+\frac{\Delta t}{2\epsilon_a}\sum\limits_{n=0}^{m-1}\llnorm{{\boldsymbol{\chi}^{n+1}+\boldsymbol{\chi}^n}}^2
		\\
		\leq&\frac{\epsilon_a}{2}\Llnorm{\dot{\boldsymbol{\theta}}}^2+\frac{\Delta t}{2\epsilon_a}4N\max_{0\leq n\leq N}\llnorm{{\boldsymbol{\chi}^n}}^2
		\\
		\leq&C\frac{\epsilon_a}{2}h^{2r}+\frac{2T}{\epsilon_a}\max_{0\leq n\leq N}\llnorm{{\boldsymbol{\chi}^n}}^2
		\end{align*}by Cauchy-Schwarz inequalities, Young's inequality and \eqref{eq:elliptic:l2} for any positive $\epsilon_a$, where $C$ depends on $\lVert\dot{\boldsymbol{w}}\rVert_{L_2(0,T;H^r(\Omega))}$.
		\item $\Delta t\sum\limits_{n=0}^{m-1}\Enorm{\boldsymbol{e}^{n+1}+\boldsymbol{e}^n}{{\boldsymbol{\chi}^{n+1}+\boldsymbol{\chi}^n}}$\\
		We follows the simple fact:
		\[\Enorm{\boldsymbol{e}^{n+1}+\boldsymbol{e}^n}{{\boldsymbol{v}}}=\Lnorm{-\nabla\cdot{\hat{\ushort{\boldsymbol{D}}}}\ushort{{\boldsymbol{\varepsilon}}}(\boldsymbol{e}^{n+1}+\boldsymbol{e}^n)}{\boldsymbol{v}}\]
		by integration by parts. Hence using Cauchy-Schwarz inequalities and Young's inequality, we can obtain
		\begin{align*}
		\Delta t&\sum_{n=0}^{m-1}\Lnorm{-\nabla\cdot{\hat{\ushort{\boldsymbol{D}}}}\ushort{{\boldsymbol{\varepsilon}}}(\boldsymbol{e}^{n+1}+\boldsymbol{e}^n)}{{\boldsymbol{\chi}^{n+1}+\boldsymbol{\chi}^n}}\\
		\leq&
		\Delta t\sum_{n=0}^{N-1}\frac{\epsilon_b}{2}\llnorm{\nabla\cdot{\hat{\ushort{\boldsymbol{D}}}}\ushort{{\boldsymbol{\varepsilon}}}(\boldsymbol{e}^{n+1}+\boldsymbol{e}^n)}^2+\frac{2T}{\epsilon_b}\max_{0\leq n\leq N}\llnorm{\boldsymbol{\chi}^n}^2
		\end{align*}for any positive $\epsilon_b$, since $m\Delta t\leq N\Delta t=T$. Recall \eqref{eq:linearinterpolation:error} then we have
		\[\llnorm{\nabla\cdot{\hat{\ushort{\boldsymbol{D}}}}\ushort{{\boldsymbol{\varepsilon}}}(\boldsymbol{e}^{n+1}+\boldsymbol{e}^n)}\leq CT^{1-\alpha}\Delta t^2,\]for some positive C depending on $\lVert\boldsymbol{w}\rVert_{C^2(0,T;H^2(\Omega))}$. Therefore, we can obtain\begin{align*}
		\Delta t\sum_{n=0}^{m-1}\Lnorm{-\nabla\cdot{\hat{\ushort{\boldsymbol{D}}}}\ushort{{\boldsymbol{\varepsilon}}}(\boldsymbol{e}^{n+1}+\boldsymbol{e}^n)}{{\boldsymbol{\chi}^{n+1}+\boldsymbol{\chi}^n}}
		\leq&CT^{3-2\alpha}\epsilon_b\Delta t^4+\frac{2T}{\epsilon_b}\max_{0\leq n\leq N}\llnorm{\boldsymbol{\chi}^n}^2.
		\end{align*}
		\item $\Delta t\sum\limits_{n=0}^{m-1}\Lnorm{\boldsymbol {\mathcal{E}}^n}{{\boldsymbol{\chi}^{n+1}+\boldsymbol{\chi}^n}}$\\
		Recalling \eqref{eq:frac:CN} and \eqref{eq:frac:CN:t0}, Cauchy-Schwarz inequalities and Young's inequality yield
		\begin{align*}
		\Delta t&\sum\limits_{n=0}^{m-1}\Lnorm{\boldsymbol {\mathcal{E}}^n}{{\boldsymbol{\chi}^{n+1}+\boldsymbol{\chi}^n}}\\
		\leq&\Delta t\sum\limits_{n=0}^{N-1}\frac{\epsilon_c}{2}\llnorm{\boldsymbol {\mathcal{E}}^n}^2+\Delta t\sum\limits_{n=0}^{N-1}\frac{2}{\epsilon_c}\max_{0\leq n\leq N}\llnorm{{\boldsymbol{\chi}^{n}}}^2\\
		\leq&\frac{\epsilon_c\Delta t^4}{8}\lVert \boldsymbol w^{(3)}\rVert_{L_2(t_1,T;L_2(\Omega))}^2+CT\epsilon_c\Delta t^{4-2\alpha}+\Delta t\sum\limits_{n=0}^{N-1}\frac{2}{\epsilon_c}\max_{0\leq n\leq N}\llnorm{{\boldsymbol{\chi}^{n}}}^2\\
		\leq&C{T\epsilon_c}\Delta t^{4-2\alpha}+\frac{2T}{\epsilon_c}\max_{0\leq n\leq N}\llnorm{{\boldsymbol{\chi}^{n}}}^2
		\end{align*}for any positive $\epsilon_c$ and some positive $C$ depending on $\lVert \boldsymbol w^{(3)}\rVert_{L_2(t_1,T;L_2(\Omega))}$ and \eqref{eq:frac:CN:t0}.
	\end{itemize}Combining the above results then \eqref{pf:frac:lemma:error:T:eq3} has a bound as
	\begin{align}
	2\rho&\llnorm{\boldsymbol{\chi}^m}^2+\frac{\Delta t^{2-\alpha}}{\Gamma(3-\alpha)}\sum_{n=0}^{m-1}\enorm{\boldsymbol{\chi}^{n+1}+\boldsymbol{\chi}^n}^2\nonumber\\
	\leq&C\left({\rho \epsilon_a}h^{2r}+T\epsilon_b \Delta t^4+ T{\epsilon_c}\Delta t^{4-2\alpha}\right)+\frac{4{\rho}T}{\epsilon_a}\max_{0\leq n\leq N}\llnorm{{\boldsymbol{\chi}^n}}^2\nonumber\\&+\frac{2T}{\epsilon_b}\max_{0\leq n\leq N}\llnorm{\boldsymbol{\chi}^n}^2+\frac{4\rho T}{\epsilon_c}\max_{0\leq n\leq N}\llnorm{{\boldsymbol{\chi}^{n}}}^2\nonumber\\&-\frac{\Delta t^{2-\alpha}}{\Gamma(3-\alpha)}\sum_{n=0}^{m-1}\Enorm{\sum_{i=0}^{n}B_{n+1,i}\boldsymbol{\chi}^i+\sum_{i=0}^{n-1}B_{n,i}\boldsymbol{\chi}^i}{\boldsymbol{\chi}^{n+1}+\boldsymbol{\chi}^n}.\label{pf:frac:lemma:error:T:eq4}
	\end{align}
	As seen in the proof of Theorem \ref{thm:frac:cg:fully:sta}, using mathematical induction we can show the bound of the last term of \eqref{pf:frac:lemma:error:T:eq4}. As proved before, coupling with $\enorm{\boldsymbol{\chi}^0}=0$, we can obtain\begin{align}
	2\rho&\llnorm{\boldsymbol{\chi}^m}^2+\frac{\Delta t^{2-\alpha}}{2\Gamma(3-\alpha)}\sum_{n=0}^{m-1}\enorm{\boldsymbol{\chi}^{n+1}+\boldsymbol{\chi}^n}^2
	\nonumber\\
	\leq&C\left({\rho \epsilon_a}h^{2r}+T^{3-2\alpha}\epsilon_b\Delta t^{4}+ T{\epsilon_c}\Delta t^{4-2\alpha}+\left(\frac{4{\rho}T}{\epsilon_a}+\frac{2T}{\epsilon_b}+\frac{4\rho T}{\epsilon_c}\right)\max_{0\leq n\leq N}\llnorm{{\boldsymbol{\chi}^{n}}}^2\right)\label{pf:frac:lemma:error:T:eq5}
	\end{align}for some positive $C$. Whence we consider maximum on \eqref{pf:frac:lemma:error:T:eq5}, we have \begin{align*}
	2\rho&\max_{0\leq n\leq N}\llnorm{\boldsymbol{\chi}^n}^2+\frac{\Delta t^{2-\alpha}}{2\Gamma(3-\alpha)}\sum_{n=0}^{N-1}\enorm{\boldsymbol{\chi}^{n+1}+\boldsymbol{\chi}^n}^2
	\nonumber\\
	\leq&2C\bigg({\rho\epsilon_a}h^{2r}+T^{3-2\alpha}\epsilon_b \Delta t^4+ T{\epsilon_c}\Delta t^{4-2\alpha}+\left(\frac{4{\rho}T}{\epsilon_a}+\frac{2T}{\epsilon_b}+\frac{4\rho T}{\epsilon_c}\right)\max_{0\leq n\leq N}\llnorm{{\boldsymbol{\chi}^{n}}}^2\bigg),
	\end{align*}therefore choosing $\epsilon_a=\epsilon_c=32CT$ and $\epsilon_b=8CT/\rho$ implies
	\begin{align*}
	\rho&\max_{0\leq n\leq N}\llnorm{\boldsymbol{\chi}^n}^2+\frac{\Delta t^{2-\alpha}}{2\Gamma(3-\alpha)}\sum_{n=0}^{N-1}\enorm{\boldsymbol{\chi}^{n+1}+\boldsymbol{\chi}^n}^2
\leq CT^{4-2\alpha}\left(h^{2r}+\Delta t^{4}+ \Delta t^{4-2\alpha}\right).
	\end{align*}As a consequence, we can conclude that
	\begin{align*}
	&\max_{0\leq n\leq N}\llnorm{\boldsymbol{\chi}^n}+\left(\Delta t^{2-\alpha}\sum_{n=0}^{N-1}\enorm{\boldsymbol{\chi}^{n+1}+\boldsymbol{\chi}^n}^2\right)^{1/2}\leq CT^{2-\alpha}(h^{r}+\Delta t^{2-\alpha}).
	\end{align*}
	
	Besides, with higher regularity of the solution in time and no singularity at $t=0$,  we could obtain second order accuracy in time. To be specific, when we suppose $\boldsymbol{w}(0)=\dot{\boldsymbol{w}}(0)=\boldsymbol{0}$ or $\boldsymbol{w}\in\ker(\mathcal{D})$, we have $H^3$ regularity in time. Therefore, instead of use of \eqref{eq:frac:CN:t0}, we can apply \eqref{eq:frac:CN:t0:opt} so that we have
	\begin{align*}
	&\max_{0\leq n\leq N}\llnorm{\boldsymbol{\chi}^n}+\left(\Delta t^{2-\alpha}\sum_{n=0}^{N-1}\enorm{\boldsymbol{\chi}^{n+1}+\boldsymbol{\chi}^n}^2\right)^{1/2}
	\leq CT^{2-\alpha}(h^{r}+\Delta t^{2}).
	\end{align*}
\end{proof}
\begin{rmk}
	Note that once the solution has higher regularity such that \[\boldsymbol{w}\in C^2(0,T;[H^{s}(\Omega)]^d)\cap W^1_\infty(0,T;\boldsymbol{V})\cap H^3(0,T;[H^{s}(\Omega)]^d),\] elliptic error estimates and \eqref{eq:frac:CN} yield
	\begin{align*}
	\max_{0\leq n\leq N}\llnorm{\boldsymbol{\chi}^n}+\left(\Delta t^{2-\alpha}\sum_{n=0}^{N-1}\enorm{\boldsymbol{\chi}^{n+1}+\boldsymbol{\chi}^n}^2\right)^{1/2}\leq CT^{2-\alpha} \lVert \boldsymbol{w}\rVert_{H^3(0,T;H^{s}(\Omega))}(h^{r}+\Delta t^{2}),
	\end{align*}for some positive $C$ depending on constants of continuity and coercivity, $\Omega$ and $\partial \Omega$ but independent of the numerical solution, mesh sizes, and time.
\end{rmk}

In the end, we can complete the error analysis as follows.
\begin{theorem}\label{thm:frac:fully:error}
	Assume that $\boldsymbol{f}$ and $\boldsymbol{w}_0$ are sufficiently smooth satisfying \textnormal{Lemma \ref{lemma:frac:error:T}}, that
	\[\boldsymbol{w}\in C^2(0,T;[H^s(\Omega)]^d)\cap W^1_\infty(0,T;\boldsymbol{V})\text{ for $s\geq2$,}\] and $\left(\boldsymbol{W}^n_h\right)_{n=0}^N$ is the fully discrete solution. Then we can observe optimal $L_2$ error as well as energy error estimates with $2-\alpha$ order accuracy in time.
	Therefore,
	\begin{align*}
	\max_{0\leq n\leq N}\llnorm{\boldsymbol{w}^n-\boldsymbol{W}_h^n}\leq 
	CT^{2-\alpha}(h^{r}+\Delta t^{2-\alpha}), \text{ and }
	\max_{0\leq n\leq N}\enorm{\boldsymbol{w}^n-\boldsymbol{W}_h^n}\leq CT^{2-\alpha}(h^{r-1}+\Delta t^{2-\alpha}),
	\end{align*}where $r=\min(k+1,s)$, for some positive $C$ independent of $h$ and $\Delta t$.
\end{theorem}
\begin{proof}
	For any $n=0,\ldots,N$, using triangular inequality, we have
	\begin{align*}
	\llnorm{\boldsymbol{w}^n-\boldsymbol{W}^n_h}=\llnorm{\boldsymbol{\theta}^n-\boldsymbol{\chi}^n}\leq\llnorm{\boldsymbol{\theta}^n}+\llnorm{\boldsymbol{\chi}^n}.
	\end{align*}
	By \eqref{eq:elliptic:l2} and Lemma \ref{lemma:frac:error:T}, it is concluded that
	\begin{align*}
	\llnorm{\boldsymbol{w}^n-\boldsymbol{W}_h^n}
	\leq& CT^{2-\alpha}(h^{r}+\Delta t^{2-\alpha}),
	\end{align*}where $C$ is from \eqref{eq:elliptic:l2} and Lemma \ref{lemma:frac:error:T}, and so on account of arbitrary $n$,
	\begin{align*}
	\max_{0\leq n\leq N}\llnorm{\boldsymbol{w}^n-\boldsymbol{W}_h^n}
	\leq& CT^{2-\alpha}(h^{r}+\Delta t^{2-\alpha}).
	\end{align*}
	In this manner, we can obtain
	\begin{align*}
	\enorm{\boldsymbol{w}^n-\boldsymbol{W}_h^n}\leq&\enorm{\boldsymbol{\theta}^n}+\enorm{\boldsymbol{\chi}^n},
	\end{align*}so that \eqref{eq:elliptic:energy} and \eqref{eq:inv} lead us to have
	\begin{align*}
	\enorm{\boldsymbol{w}^n-\boldsymbol{W}_h^n}\leq&\enorm{\boldsymbol{\theta}^n}+Ch^{-1}\llnorm{\boldsymbol{\chi}^n}\leq CT^{2-\alpha}(h^{r-1}+\Delta t^{2-\alpha}).
	\end{align*}
	Thus, we have
	\[  \max_{0\leq n\leq N}\enorm{\boldsymbol{w}^n-\boldsymbol{W}^n_h}\leq CT^{2-\alpha}(h^{r-1}+\Delta t^{2-\alpha}).\]
\end{proof}
\begin{corollary}\label{cor:frac:cg:error}
	Under the same conditions in \textnormal{Theorem \ref{thm:frac:fully:error}}, we suppose higher regularity in time such that $\boldsymbol{{w}}\in H^3(0,T;[H^s(\Omega)]^d)$ or we further assume that \eqref{eq:frac:CN:t0:opt} is satisfied. Then we can obtain optimal results of Crank-Nicolson scheme i.e.,\begin{align*}
	\max_{0\leq n\leq N}\llnorm{\boldsymbol{w}^n-\boldsymbol{W}_h^n}\leq CT^{2-\alpha}\lVert \boldsymbol{w}\rVert_{H^3(0,T;H^{r}(\Omega))}(h^{r}+\Delta t^{2})&,\\
	\max_{0\leq n\leq N}\enorm{\boldsymbol{w}^n-\boldsymbol{W}_h^n}\leq CT^{2-\alpha}\lVert \boldsymbol{w}\rVert_{H^3(0,T;H^{r}(\Omega))}(h^{r-1}+\Delta t^{2})&,
	\end{align*}where $C$ is a positive constant such that is independent of solutions, mesh sizes, $T$ but depends on the domain, its boundary and coefficients of coercivity and continuity.
\end{corollary}
\begin{proof}
	As shown in Theorem \ref{thm:frac:fully:error}, triangular inequalities combined with \eqref{eq:elliptic:energy}, \eqref{eq:elliptic:l2} and Lemma \ref{lemma:frac:error:T} complete the proof.
\end{proof}
\section{Numerical Experiments}\label{sec:numeric}
\noindent We have carried out numerical experiments using FEniCS (\url{https://fenicsproject.org/}). In this section, we present tables of numerical errors, as well as convergence rates for some evidence of the above error estimates theorem in practice. Codes are available at the author's Github (\url{https://github.com/Yongseok7717/Visco_Frac_CG}) written as python scripts to reproduce the tabulated results and figures that are given below. In addition, using Docker container, we can also run the codes at a \texttt{bash prompt}, e.g. the commands to run are
	\begin{center}
		{\texttt{docker pull variationalform/fem:yjcg2\\
		docker run -ti variationalform/fem:yjcg2\\
		cd; cd  .\/codes\/Visco\_Frac\_CG-master; .\/main.sh}}
\end{center}
Consider two cases;  one is an example that is not of class $H^3$ in time but the other is a smoother case. We set our spatial domain as the unit square, $T=1$ and $\alpha=1/2$.

\begin{figure}[H]
	\centering
		\includegraphics[width=0.55\textwidth]{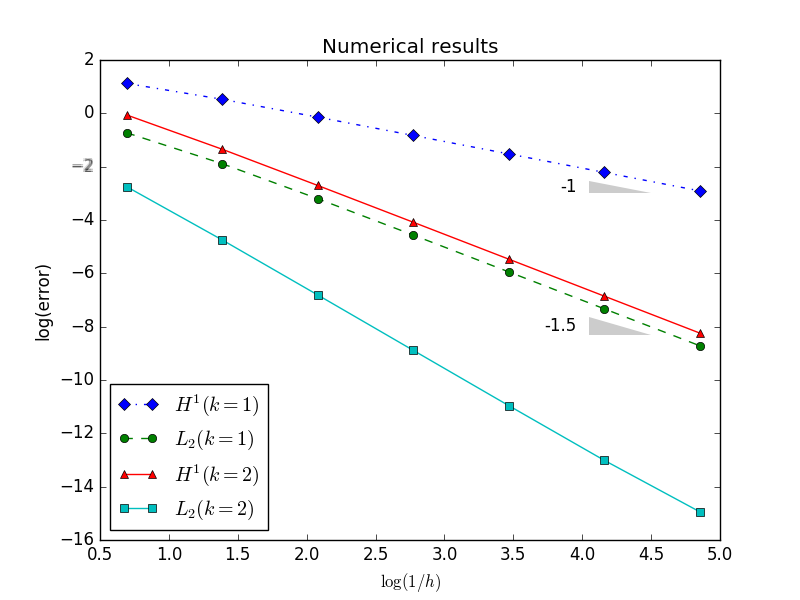}
	\caption{Numerical convergent orders of Example \ref{eg:nonsmooth}: \textbf{linear} (dash line) and \textbf{quadratic} (solid line) polynomial basis}
	\label{fig:frac:nonsmooth}
\end{figure}

\begin{figure}[H]\vspace*{-0.5cm}
	\centering
	\includegraphics[width=0.55\textwidth]{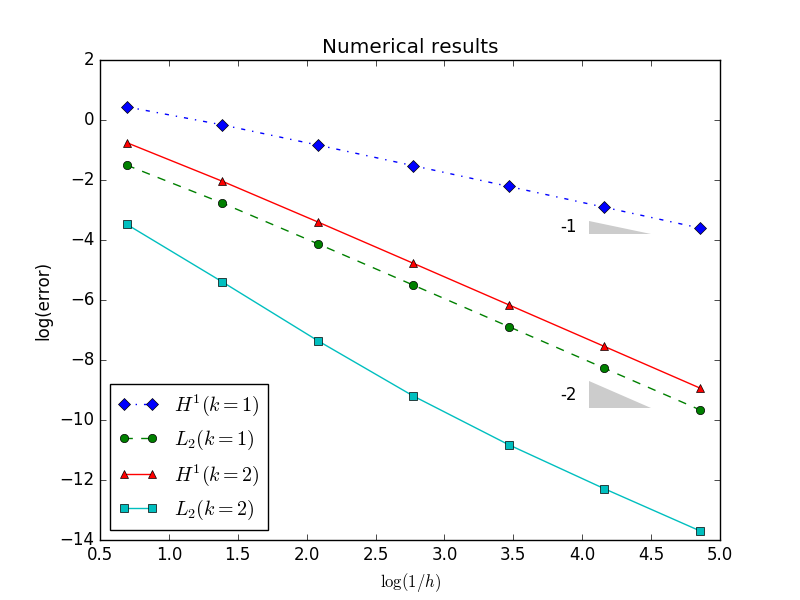}
	\caption{Numerical convergent orders of Example \ref{eg:smooth}: \textbf{linear} (dash line) and \textbf{quadratic} (solid line) polynomial basis}
	\label{fig:frac:smooth}
\end{figure}

\begin{eg}\label{eg:nonsmooth}	
	Let us define \begin{align*}\boldsymbol{w}(x,y,t)=(t+t^{1.5})\left[\begin{array}{c}
	\sin(\pi x)\sin(\pi y)  \\
	xy(1-x)(1-y)
	\end{array}
	\right].
	\end{align*}Then $\boldsymbol{w}\in C^2(0,T;[C^\infty(\Omega)]^2)\cap W^2_1(0,T;[C^\infty(\Omega)]^2)$ with homogeneous Dirichlet boundary condition. Also, we can derive data terms which satisfy \eqref{eq:frac:primal:eq1}. Note that $\boldsymbol{w}^{(3)}(t)$ is not bounded and not integrable in time so that we cannot fully take an advantage of second order schemes. However, we can observe suboptimal results but higher than first order schemes.
	
	Let us define $\boldsymbol{e}^n=\boldsymbol{w}(t_n)-\boldsymbol{W}_h^n$ for $n=0,\ldots,N$. 
	By error estimates theorems for both solutions, we have
	\[\enorm{\boldsymbol{e}^n}=O(h^{k}+\Delta t^{1.5}),\qquad\text{and}\qquad\llnorm{\boldsymbol{e}^n}=O(h^{k+1}+\Delta t^{1.5}),\]since $s=\infty$. In other words, the orders of convergence depend only on the degree of polynomial $k$ for the spatial mesh. On the other hand, regardless of types of the norm, convergence rates of time are suboptimally fixed by 1.5.
	
	In Tables \ref{table:vector:frac:nonsmooth:linear} and \ref{table:vector:frac:nonsmooth:quadratic}, we can observe $H^1$ norm and $L_2$ norm errors for linear and quadratic polynomial bases, respectively. Also, we can observe the numerical convergent order with respect to the polynomial degrees of $k$ when $\Delta t$ is sufficiently small in Table \ref{table:convrate:nonsmooth}. In a similar way, we could compute the rate of convergence with respect to time for small $h$. However, in a practical sense, it is difficult to computationally solve it for fine meshes if the machine is not sufficiently good enough. In other words, we may encounter some memory issues. For example, when $h=1/512$, $\hnorm{\boldsymbol{e}}^N$ are given by 1.442e-2 and 1.368e-2, for $\Delta t=1/4$ and $\Delta t=1/8$, respectively. It implies that $h=1/512$ is not small enough to see the convergent order of time but smaller spatial meshes enforce us to have large systems of matrix and memory issues. Alternatively, while we consider $\Delta t\approx h$, the numerical convergent rate $d_c$ can be computed by
	$d_c\ =\ [{\log(\text{error of }h_1)-\log(\text{error of }h_2)}]/[{\log(h_1)-\log(h_2)}].$
	Here, $d_c$ can represent the convergent order of time if we consider $k\geq2$ or $L_2$ norm errors. For example, when we take diagonals of Tables \ref{table:vector:frac:nonsmooth:linear} and \ref{table:vector:frac:nonsmooth:quadratic}, the convergent rates are illustrated as the gradients of line in Figure \ref{fig:frac:nonsmooth}. For the linear polynomial basis, the numerical rate of the energy norm (equivalent to $H^1$ norm) is $d_c\approx1$, otherwise $d_c\approx 1.5$ for higher degree of polynomial or $L_2$ norm. 
	Interestingly, in Figure \ref{fig:frac:nonsmooth}, the slope of line for $L_2$ errors of $k=2$ looks steeper than 1.5. Theoretically, we can rewrite the $L_2$ norm of error for this case as 
	$\llnorm{\boldsymbol{e}^n}\approx C_1 h^3+C_2h^{1.5}$. Hence if $h$ is not small enough, $d_c$ could be greater than 1.5. However, as $h$ decreasing, $d_c$ will approach to 1.5. For instance, the setting of $h=\Delta t$ and $k=2$ leads us to obtain Table \ref{table:L2quad:nonsmooth} which exhibits that the convergent orders of $\llnorm{\boldsymbol{e}^N}$ is higher than 1.5 but they are decreasing while $h$ becomes smaller.
\end{eg} \vspace*{-0.5cm}
\begin{table}[H]
	\footnotesize
	\centering
	$H^1$ error \linebreak
	\begin{tabular}{|c|ccccccc|}\hline
		\backslashbox{$h$\kern-1em}{\kern-1em$\Delta t$}	&          1/8  &         1/16  &         1/32  &         1/64  &        1/128  &        1/256  &        1/512   \\\hline
		1/2 &   3.072& 3.072& 3.073& 3.073& 3.073& 3.073& 3.073\\
		1/4& 1.694 &1.694 &1.694& 1.694 &1.694& 1.694 &1.694\\
		1/8& 8.677e-01& 8.677e-01 &8.677e-01& 8.677e-01& 8.677e-01& 8.677e-01& 8.677e-01\\
		1/16& 4.364e-01& 4.364e-01& 4.364e-01& 4.364e-01& 4.364e-01& 4.364e-01& 4.364e-01\\
		1/32& 2.185e-01& 2.185e-01& 2.185e-01& 2.185e-01& 2.185e-01& 2.185e-01& 2.185e-01\\
		1/64& 1.093e-01& 1.093e-01 &1.093e-01 &1.093e-01 &1.093e-01 &1.093e-01 &1.093e-01\\
		1/128& 5.466e-02& 5.465e-02& 5.465e-02 &5.465e-02 &5.465e-02& 5.465e-02& 5.465e-02  \\
		\hline
	\end{tabular}\vspace{0.1cm}

$L_2$ error \linebreak
\begin{tabular}{|c|ccccccc|}\hline
	\backslashbox{$h$\kern-1em}{\kern-1em$\Delta t$}&          1/8  &         1/16  &         1/32  &         1/64  &        1/128  &        1/256  &        1/512   \\\hline
	1/2 &  4.827e-01& 4.824e-01 &4.823e-01& 4.823e-01& 4.823e-01 &4.823e-01& 4.823e-01\\
	1/4 &1.519e-01 &1.515e-01 &1.513e-01 &1.513e-01 &1.513e-01 &1.513e-01& 1.513e-01\\
	1/8 &4.103e-02 &4.087e-02 &4.080e-02 &4.079e-02 &4.078e-02 &4.078e-02& 4.078e-02\\
	1/16 &1.057e-02 &1.049e-02 &1.045e-02 &1.043e-02 &1.043e-02 &1.043e-02& 1.043e-02\\
	1/32 &2.744e-03 &2.679e-03 &2.640e-03 &2.627e-03 &2.624e-03 &2.622e-03& 2.622e-03\\
	1/64 &7.745e-04 &7.125e-04 &6.738e-04 &6.616e-04 &6.581e-04 &6.570e-04& 6.567e-04\\
	1/128 &2.864e-04 &2.215e-04 &1.815e-04 &1.693e-04 &1.657e-04 &1.647e-04& 1.643e-04\\
		\hline\end{tabular}
	\caption{Numerical errors; Example \ref{eg:nonsmooth}; $k=1,\ n=N$}
	\label{table:vector:frac:nonsmooth:linear}
\end{table}\vspace{-1.0cm}
\begin{table}[H]
	\footnotesize
	\centering
	$H^1$ error \linebreak\begin{tabular}{|c|ccccccc|}\hline
		\backslashbox{$h$\kern-1em}{\kern-1em$\Delta t$}& 1/8  & 1/16  & 1/32  & 1/64  & 1/128  & 1/256&1/512   \\\hline
		1/2 &   9.417e-01& 9.417e-01& 9.417e-01& 9.417e-01& 9.417e-01& 9.417e-01& 9.417e-01\\
		1/4 &2.604e-01 &2.604e-01 &2.604e-01& 2.604e-01& 2.604e-01& 2.604e-01& 2.604e-01\\
		1/8& 6.700e-02 &6.700e-02 &6.700e-02 &6.700e-02 &6.700e-02 &6.700e-02& 6.700e-02\\
		1/16& 1.689e-02 &1.688e-02 &1.688e-02 &1.688e-02 &1.688e-02 &1.688e-02& 1.688e-02\\
		1/32& 4.276e-03 &4.238e-03 &4.229e-03 &4.228e-03 &4.228e-03& 4.228e-03& 4.228e-03\\
		1/64& 1.236e-03 &1.098e-03 &1.061e-03 &1.058e-03 &1.058e-03 &1.057e-03& 1.057e-03\\
		1/128& 6.922e-04& 3.954e-04 &2.796e-04 &2.658e-04 &2.645e-04 &2.644e-04& 2.644e-04\\
		\hline\end{tabular}\vspace{0.1cm}
	
	$L_2$ error \linebreak\begin{tabular}{|c|ccccccc|}\hline
		\backslashbox{$h$\kern-1em}{\kern-1em$\Delta t$}&1/8  &1/16  &1/32  &1/64  &1/128  &1/256  &1/512 \\\hline
	1/2 &  6.394e-02& 6.382e-02 &6.377e-02 &6.375e-02 &6.375e-02 &6.375e-02 &6.375e-02\\
	1/4& 8.718e-03 &8.688e-03 &8.671e-03 &8.665e-03 &8.664e-03 &8.663e-03 &8.663e-03\\
	1/8& 1.133e-03 &1.114e-03 &1.104e-03 &1.101e-03 &1.100e-03 &1.100e-03 &1.100e-03\\
	1/16& 2.013e-04 &1.577e-04 &1.412e-04 &1.385e-04 &1.380e-04 &1.379e-04 &1.378e-04\\
	1/32& 1.358e-04 &6.726e-05 &2.699e-05 &1.851e-05 &1.742e-05 &1.727e-05 &1.724e-05\\
	1/64 &1.339e-04 &6.424e-05 &2.008e-05 &6.376e-06 &2.844e-06 &2.240e-06 &2.166e-06\\
	1/128 &1.338e-04& 6.416e-05 &1.992e-05 &5.956e-06 &1.826e-06 &6.248e-07 &3.251e-07\\
		\hline\end{tabular}
	\caption{Numerical errors; Example \ref{eg:nonsmooth}; $k=2,\ n=N$}
	\label{table:vector:frac:nonsmooth:quadratic}
\end{table}
\vspace*{-1.0cm}
\begin{table}[H]
	\footnotesize
	\centering
	\begin{tabular}{|c|cccc|cccc|}\hline
		\multirow{2}{*}{$h$}	&          \multicolumn{4}{c|}{$k=1$}&\multicolumn{4}{c|}{$k=1$}\\
		&$\hnorm{\boldsymbol{e}^N}$&Rate&$\llnorm{\boldsymbol{e}^N}$&Rate&$\hnorm{\boldsymbol{e}^N}$&Rate&$\llnorm{\boldsymbol{e}^N}$&Rate\\
		\hline
		1/2&3.073&&4.823e-01&&9.417e-01&&6.375e-02&\\
		1/4&1.694&0.86&1.513e-01&1.67&2.604e-01&1.85&8.663e-03&2.88\\
		1/8&8.677e-01&0.97&4.078e-02&1.89&6.700e-02&1.96&1.100e-03&2.98\\
		1/16&4.364e-01&0.99&1.043e-02&1.97&1.688e-02&1.99&1.378e-04&3.00\\
		1/32&2.185e-01&1.00&2.622e-03&1.99&4.228e-03&2.00&1.724e-05&3.00\\
		\hline
	\end{tabular}
	\caption{Convergent rates; Example \ref{eg:nonsmooth}; $\Delta t=1/512$}
	\label{table:convrate:nonsmooth}
\end{table}
\vspace*{-1cm}
\vspace*{-0.5cm}\begin{table}[H]
	\footnotesize
	\centering
	\begin{tabular}{|c|cccccc|}\hline
	$h$&1/8&1/16&1/32&1/64&1/128&1/256\\
		\hline
		Error(rate)&1.133e-03&1.577e-04(2.85)&2.699e-05(2.55)&6.376e-06(2.08)&1.826e-06(1.80)&5.616e-07(1.70)\\
		\hline
	\end{tabular}
	\caption{$L_2$ norm errors of Example \ref{eg:nonsmooth} for $\Delta t= h$ and $k=2$}
	\label{table:L2quad:nonsmooth}
\end{table}\vspace*{-0.0cm}
Due to loss of $H^3$ regularity in time, Example \ref{eg:nonsmooth} cannot take fully the advantage of second order scheme. However, once we give further assumptions for higher regularity such as $\boldsymbol{w}(0)=\dot{\boldsymbol{w}}(0)=0$, our fully discrete formulation will guarantee spatially optimal error estimates as well as second order accuracy in time.
\begin{eg}\label{eg:smooth}	
	Let  \begin{align*}\boldsymbol{w}(x,y,t)=t^{3.5}\left[\begin{array}{c}
	\sin(\pi x)\sin(\pi y)  \\
	xy(1-x)(1-y)
	\end{array}
	\right].
	\end{align*}
	The exact solution is of class $C^3$ in time, i.e. Example \ref{eg:smooth} has higher regularity than Example \ref{eg:nonsmooth} with respect to time. Therefore, according to Corollary \ref{cor:frac:cg:error}, we have
	\[\enorm{\boldsymbol{e}^n}=O(h^{k}+\Delta t^{2}),\qquad\text{and}\qquad\llnorm{\boldsymbol{e}^n}=O(h^{k+1}+\Delta t^{2}).\]
	
	Tables \ref{table:vector:frac:smooth:linear} and \ref{table:vector:frac:smooth:quadratic} indicate that the orders of spatial convergence are optimal not only in $H^1$ norm but also in $L_2$ norm. More precisely, we can observe the convergent order in Table \ref{table:convrate:smooth}.  Furthermore, when $\Delta t\approx h$, we can observe numerical convergent rates in Figure \ref{fig:frac:smooth}. The energy error estimates show first order for the linear polynomial basis. On the other hand, regardless of a degree of polynomials, $L_2$ norm errors have second order accuracy, i.e. $d_c\approx2$. 
\end{eg}

Comparing Example \ref{eg:nonsmooth} and Example \ref{eg:smooth}, we can observe optimal error estimates with respect to space but not enough regularity in time restricts the convergence order of time. Nevertheless, sufficiently smooth data terms enable our numerical scheme to have better accuracy than first order finite difference methods, e.g. it is of order $2-\alpha$. In addition, once we assume $H^3$ regularity in time, we get the second order of convergence in time.

\begin{table}[H]
	\footnotesize
	\centering
	$H^1$ error \linebreak
		\begin{tabular}{|c|ccccccc|}\hline
		\backslashbox{$h$}{$\Delta t$}	         & 1/8  & 1/16  & 1/32  & 1/64  & 1/128  & 1/256  & 1/512   \\\hline
		1/2 &  1.551 &1.553 &1.553 &1.553 &1.553 &1.553 &1.553\\
		1/4& 8.492e-01& 8.505e-01& 8.509e-01 &8.510e-01 &8.510e-01 &8.510e-01& 8.510e-01\\
		1/8& 4.337e-01& 4.341e-01& 4.343e-01 &4.343e-01 &4.343e-01 &4.344e-01& 4.344e-01\\
		1/16& 2.185e-01& 2.182e-01& 2.182e-01 &2.183e-01 &2.183e-01 &2.183e-01& 2.183e-01\\
		1/32& 1.105e-01& 1.093e-01& 1.093e-01 &1.093e-01 &1.093e-01 &1.093e-01& 1.093e-01\\
		1/64& 5.755e-02& 5.483e-02 &5.465e-02 &5.465e-02 &5.465e-02 &5.465e-02& 5.465e-02\\
		1/128& 3.291e-02& 2.773e-02 &2.735e-02 &2.733e-02 &2.733e-02 &2.733e-02& 2.733e-02  \\
				\hline
	\end{tabular}\vspace{0.1cm}
	
	$L_2$ error\linebreak
		\begin{tabular}{|c|ccccccc|}\hline
		\backslashbox{$h$}{$\Delta t$}&          1/8  &1/16  &1/32  &1/64  &1/128  &1/256  &1/512   \\\hline
		1/2 &2.224e-01& 2.212e-01& 2.208e-01& 2.208e-01& 2.207e-01 &2.207e-01& 2.207e-01\\
		1/4& 6.496e-02& 6.286e-02& 6.232e-02& 6.218e-02& 6.214e-02 &6.213e-02& 6.213e-02\\
		1/8& 1.923e-02& 1.681e-02& 1.620e-02& 1.604e-02 &1.600e-02 &1.599e-02& 1.599e-02\\
		1/16& 7.605e-03& 4.898e-03& 4.245e-03& 4.083e-03 &4.043e-03 &4.033e-03& 4.030e-03\\
		1/32& 4.875e-03& 1.953e-03& 1.235e-03& 1.065e-03 &1.023e-03 &1.013e-03& 1.010e-03\\
		1/64& 4.249e-03& 1.268e-03& 4.963e-04& 3.101e-04 &2.665e-04 &2.560e-04& 2.534e-04\\
		1/128& 4.099e-03& 1.111e-03& 3.250e-04& 1.254e-04 &7.774e-05 &6.668e-05& 6.401e-05\\
		\hline\end{tabular}
	\caption{Numerical errors; Example \ref{eg:smooth}; $k=1,\ n=N$}
	\label{table:vector:frac:smooth:linear}
\end{table}\vspace*{-1cm}
\begin{table}[H]
	\footnotesize
	\centering
	$H^1$ error \linebreak\begin{tabular}{|c|ccccccc|}\hline
		\backslashbox{$h$}{$\Delta t$}& 1/8  & 1/16  & 1/32  & 1/64  & 1/128  & 1/256  & 1/512   \\\hline
		1/2 &  4.707e-01 &4.712e-01& 4.714e-01 &4.715e-01& 4.715e-01& 4.715e-01& 4.715e-01\\
		1/4 &1.312e-01 &1.302e-01& 1.302e-01 &1.302e-01 &1.302e-01& 1.302e-01& 1.302e-01\\
		1/8 &3.817e-02 &3.383e-02& 3.352e-02 &3.350e-02 &3.350e-02& 3.350e-02& 3.350e-02\\
		1/16 &2.028e-02 &9.720e-03 &8.529e-03 &8.445e-03 &8.439e-03& 8.439e-03& 8.439e-03\\
		1/32 &1.857e-02 &5.275e-03 &2.453e-03 &2.138e-03 &2.115e-03& 2.114e-03& 2.114e-03\\
		1/64 &1.846e-02 &4.862e-03 &1.353e-03& 6.168e-04 &5.348e-04& 5.291e-04& 5.288e-04\\
		1/128 &1.845e-02& 4.835e-03 &1.252e-03& 3.441e-04 &1.548e-04& 1.338e-04& 1.323e-04\\		
		\hline\end{tabular}\vspace{0.1cm}
	
	$L_2$ error \linebreak\begin{tabular}{|c|ccccccc|}\hline
		\backslashbox{$h$}{$\Delta t$} &1/8  &1/16  &1/32  &1/64  &1/128  &1/256  &1/512   \\\hline
		1/2 &   3.095e-02& 2.940e-02& 2.902e-02 &2.893e-02 &2.891e-02 &2.890e-02& 2.890e-02\\
		1/4& 6.520e-03& 4.552e-03& 4.236e-03& 4.174e-03& 4.160e-03& 4.157e-03& 4.156e-03\\
		1/8& 4.160e-03& 1.257e-03& 6.408e-04& 5.569e-04& 5.456e-04& 5.435e-04& 5.430e-04\\
		1/16& 4.055e-03& 1.068e-03& 2.865e-04& 1.013e-04& 7.210e-05& 6.912e-05& 6.875e-05\\
		1/32& 4.050e-03& 1.061e-03& 2.738e-04& 7.056e-05& 1.994e-05& 9.837e-06& 8.721e-06\\
		1/64& 4.050e-03& 1.061e-03& 2.733e-04& 6.976e-05& 1.773e-05& 4.610e-06& 1.570e-06\\
		1/128& 4.050e-03& 1.061e-03& 2.733e-04& 6.973e-05& 1.768e-05& 4.466e-06& 1.133e-06
		\\		
		\hline\end{tabular}
	\caption{Numerical errors; Example \ref{eg:smooth}; $k=2\ n=N$}
	\label{table:vector:frac:smooth:quadratic}
\end{table}
\vspace*{-1cm}
\begin{table}[H]
\footnotesize
\centering
\begin{tabular}{|c|cccc|cccc|}\hline
	\multirow{2}{*}{$h$}	&          \multicolumn{4}{c|}{$k=1$}&\multicolumn{4}{c|}{$k=1$}\\
	&$\hnorm{\boldsymbol{e}^N}$&Rate&$\llnorm{\boldsymbol{e}^N}$&Rate&$\hnorm{\boldsymbol{e}^N}$&Rate&$\llnorm{\boldsymbol{e}^N}$&Rate\\
	\hline
	1/2&1.553&&2.207e-01&&4.715e-01&&2.890e-02&\\
	1/4&8.510e-01&0.87&6.213e-02&1.83&1.302e-01&1.86&4.156e-03&2.80\\
	1/8&4.344e-01&0.97&1.599e-02&1.96&3.350e-02&1.96&5.430e-04&2.94\\
	1/16&2.183e-01&0.99&4.030e-03&1.99&8.439e-03&1.99&6.875e-05&2.98\\
	1/32&1.093e-01&1.00&1.010e-03&2.00&2.114e-03&2.00&8.721e-06&2.98\\
	\hline
\end{tabular}
\caption{Convergent rates; Example \ref{eg:nonsmooth}; $\Delta t=1/512$}
\label{table:convrate:smooth}
\end{table}
\section{Conclusion}\label{sec:conclusion}
In conclusion, the numerical scheme of the fractional order viscoelasticity problem has been formulated. Without Gr\"onwall's inequality, we can show stability bounds for semi-discrete and fully discrete schemes which are non-exponentially increasing with respect to the final time. \textit{A priori} error estimates have been derived for the fully discrete formulation. We gives a remark regarding regularity of solution in time, which restricts $H^3$ smoothness in time due to weak singularity. However, we can take some advantage of second order schemes in time where we assume smooth data, and higher regularity enables the order of convergence optimal in time. In the end, we have illustrated numerical examples of suboptimal and optimal cases.
\bibliography{ref}

\begin{thebibliography}{10}

\bibitem{hunter1976mechanics}
S.~C. Hunter, {\em Mechanics of continuous media}.
\newblock Halsted Press, 1976.

\bibitem{VE}
S.~Shaw and J.~Whiteman, ``{Some partial differential Volterra equation
  problems arising in viscoelasticity},'' in {\em Proceedings of Equadiff},
  vol.~9, pp.~183--200, 1998.

\bibitem{DGV}
B.~Rivi{\`e}re, S.~Shaw, and J.~Whiteman, ``{Discontinuous {G}alerkin finite
  element methods for dynamic linear solid viscoelasticity problems},'' {\em
  Numerical Methods for Partial Differential Equations}, vol.~23, no.~5,
  pp.~1149--1166, 2007.

\bibitem{findley2013creep}
W.~N. Findley and F.~A. Davis, {\em Creep and relaxation of nonlinear
  viscoelastic materials}.
\newblock Courier Corporation, 2013.

\bibitem{drozdov1998viscoelastic}
A.~D. Drozdov, {\em Viscoelastic structures: mechanics of growth and aging}.
\newblock Academic Press, 1998.

\bibitem{golden2013boundary}
J.~M. Golden and G.~A. Graham, {\em Boundary value problems in linear
  viscoelasticity}.
\newblock Springer Science \& Business Media, 2013.

\bibitem{riviere2003discontinuous}
B.~Rivi\'ere, S.~Shaw, M.~F. Wheeler, and J.~R. Whiteman, ``{Discontinuous
  Galerkin finite element methods for linear elasticity and quasistatic linear
  viscoelasticity},'' {\em Numerische Mathematik}, vol.~95, no.~2,
  pp.~347--376, 2003.

\bibitem{shaw1998numerical}
S.~Shaw and J.~Whiteman, ``Numerical solution of linear quasistatic hereditary
  viscoelasticity problems ii: a posteriori estimates,'' tech. rep., BICOM
  Technical Report 98-3, see www. brunel. ac. uk/\~{} icsrbicm, 1998.

\bibitem{shaw1999numerical}
S.~Shaw and J.~Whiteman, ``Numerical solution of linear quasistatic hereditary
  viscoelasticity problems i: a priori estimates,'' {\em Recall}, vol.~11,
  p.~4, 1999.

\bibitem{jang2020finite}
Y.~Jang and S.~Shaw, ``Finite element approximation and analysis of
  viscoelastic wave propagation with internal variable formulations,'' {\em
  arXiv preprint arXiv:2001.04745}, 2020.

\bibitem{nutting1921new}
P.~Nutting, ``A new general law of deformation,'' {\em Journal of the Franklin
  Institute}, vol.~191, no.~5, pp.~679--685, 1921.

\bibitem{torvik1984appearance}
P.~J. Torvik and R.~L. Bagley, ``On the appearance of the fractional derivative
  in the behavior of real materials,'' {\em Journal of Applied Mechanics},
  vol.~51, no.~2, pp.~294--298, 1984.

\bibitem{koeller1984applications}
R.~Koeller, ``Applications of fractional calculus to the theory of
  viscoelasticity,'' vol.~51, no.~2, pp.~299--307, 1984.

\bibitem{li2011numerical}
C.~Li, A.~Chen, and J.~Ye, ``Numerical approaches to fractional calculus and
  fractional ordinary differential equation,'' {\em Journal of Computational
  Physics}, vol.~230, no.~9, pp.~3352--3368, 2011.

\bibitem{mclean1993numerical}
W.~McLean and V.~Thom{\'e}e, ``Numerical solution of an evolution equation with
  a positive-type memory term,'' {\em The ANZIAM Journal}, vol.~35, no.~1,
  pp.~23--70, 1993.

\bibitem{mclean2010numerical}
W.~McLean and V.~Thom{\'e}e, ``{Numerical solution via Laplace transforms of a
  fractional order evolution equation},'' {\em The Journal of Integral
  Equations and Applications}, pp.~57--94, 2010.

\bibitem{mclean2010maximum}
W.~McLean and V.~Thom{\'e}e, ``Maximum-norm error analysis of a numerical
  solution via laplace transformation and quadrature of a fractional-order
  evolution equation,'' {\em IMA journal of numerical analysis}, vol.~30,
  no.~1, pp.~208--230, 2010.

\bibitem{linz2001theoretical}
P.~Linz, {\em Theoretical numerical analysis: an introduction to advanced
  techniques}.
\newblock Courier Corporation, 2001.

\bibitem{oldham1974fractional}
K.~Oldham and J.~Spanier, {\em The fractional calculus theory and applications
  of differentiation and integration to arbitrary order}, vol.~111.
\newblock Elsevier, 1974.

\bibitem{malinowska2012introduction}
A.~B. Malinowska and D.~F. Torres, {\em Introduction to the fractional calculus
  of variations}.
\newblock World Scientific Publishing Company, 2012.

\bibitem{miller1993introduction}
K.~S. Miller and B.~Ross, {\em An introduction to the fractional calculus and
  fractional differential equations}.
\newblock Wiley-Interscience, 1993.

\bibitem{MTE}
S.~Brenner and R.~Scott, {\em The mathematical theory of finite element
  methods}, vol.~15.
\newblock Springer Science \& Business Media, 2007.

\bibitem{wheeler}
M.~F. Wheeler, ``{A priori ${L}_2$ error estimates for Galerkin approximations
  to parabolic partial differential equations},'' {\em SIAM Journal on
  Numerical Analysis}, vol.~10, no.~4, pp.~723--759, 1973.

\bibitem{ciarlet2010korn}
P.~G. Ciarlet, ``{On Korn’s inequality},'' {\em Chinese Annals of
  Mathematics, Series B}, vol.~31, no.~5, pp.~607--618, 2010.

\bibitem{horgan1983inequalities}
C.~O. Horgan and L.~E. Payne, ``{On inequalities of Korn, Friedrichs and
  Babu{\v{s}}ka-Aziz},'' {\em Archive for Rational Mechanics and Analysis},
  vol.~82, no.~2, pp.~165--179, 1983.

\bibitem{nitsche1981korn}
J.~A. Nitsche, ``{On Korn's second inequality},'' {\em RAIRO. Analyse
  num{\'e}rique}, vol.~15, no.~3, pp.~237--248, 1981.

\bibitem{li2011developing}
J.~Li, Y.~Huang, and Y.~Lin, ``{Developing finite element methods for Maxwell's
  equations in a Cole--Cole dispersive medium},'' {\em SIAM Journal on
  scientific computing}, vol.~33, no.~6, pp.~3153--3174, 2011.

\bibitem{WARBURTON20032765}
T.~Warburton and J.~Hesthaven, ``On the constants in hp-finite element trace
  inverse inequalities,'' {\em Computer Methods in Applied Mechanics and
  Engineering}, vol.~192, no.~25, pp.~2765 -- 2773, 2003.

\bibitem{DG}
B.~Rivi{\`e}re, {\em {Discontinuous Galerkin methods for solving elliptic and
  parabolic equations: theory and implementation}}.
\newblock SIAM, 2008.

\bibitem{ozisik2010constants}
S.~Ozisik, B.~Riviere, and T.~Warburton, ``{On the constants in inverse
  inequalities in $L_2$},'' tech. rep., Rice University, 2010.

\bibitem{thomee1984galerkin}
V.~Thom{\'e}e, {\em Galerkin finite element methods for parabolic problems},
  vol.~1054.
\newblock Springer, 1984.

\bibitem{grenander1958toeplitz}
U.~Grenander and G.~Szeg{\"o}, {\em Toeplitz forms and their applications}.
\newblock Univ of California Press, 1958.

\bibitem{lopez1990difference}
J.~Lopez-Marcos, ``A difference scheme for a nonlinear partial
  integrodifferential equation,'' {\em SIAM journal on numerical analysis},
  vol.~27, no.~1, pp.~20--31, 1990.

\bibitem{dauge1988elliptic}
M.~Dauge, ``Elliptic boundary value problems on corner domains, volume 1341 of
  lecture notes in mathematics,'' 1988.

\bibitem{grisvard2011elliptic}
P.~Grisvard, {\em Elliptic problems in nonsmooth domains}.
\newblock SIAM, 2011.

\end{thebibliography}
\bibliographystyle{ieeetr}
\end{document}